\documentclass[11pt]{article}
\usepackage{amssymb,amsfonts,amsmath,amsthm,mathrsfs,enumerate,cancel,graphicx}
\usepackage{varwidth, float,enumitem,makecell,caption,color}
 \usepackage{tikz}
 \usepackage{pgflibrarysnakes}

\usepackage{graphicx}
\usepackage{color}
\usepackage[colorlinks = true,
linkcolor = blue,
urlcolor  = blue,
citecolor = blue,
anchorcolor = blue]{hyperref}

\numberwithin{equation}{section}

\begin{document}

\title{\bf
Fundamentals of thermoelasticity for curved beams	
%
%
%
%
%
%
}

\author{
	\small {\bf M. A. Jorge Silva\thanks{Supported by the CNPq, grant  160954/2022-3. Email: 
			\href{mailto:marcioajs@uel.br}{marcioajs@uel.br}.  (Corresponding Author)}}\\
	\small Department of Mathematics, State University of Londrina, \\
	\small 86057-970 Londrina, PR, Brazil. \\
	\small {\bf T. F. Ma\thanks{Supported by the CNPq grant 315165/2021-9 and FAPDF grant 193.00001821/2022-21. \\   Email: 
			\href{mailto:matofu@mat.unb.br}{matofu@mat.unb.br}.}}\\
	\small Department of Mathematics, University of Bras\'ilia,\\
	\small 70910-900 Bras\'ilia, DF, Brazil. \\
}
\date{}
\maketitle

\begin{abstract} 
The purpose of this paper is twofold. Firstly, we conduct an in-depth analysis of mathematical modeling concerning thermal-mechanical curved beams, by taking into consideration three primary forces widely accepted in the literature: axial load, shear force, and bending moment. Additionally, we examine their appropriate thermal couplings, shedding light on the intricate interplay between stress-strain relationships and temperature variations. This analysis is situated within the well-recognized context of the Bresse governing model for arched beams.
Secondly, drawing upon distinguished constitutive laws for heat flux of conduction, we compile a comprehensive list of thermoelastic curved beam systems in various scenarios. We introduce new categories of problems that exhibit specific features from the  thermal  point of view.
\end{abstract}

\noindent{\bf Keywords:} Curved beams, thermoelasticity, constitutive law.

\noindent{\bf 2020 MSC:} 35Q79, 74A15, 74F05,   74K10, 80A05.

%

\section{Introduction}
\label{sec-introd}

This work is mainly focused with mathematical models for arched beams. It is known 
that arches and similar curved materials have been utilized since ancient civilizations, including the Assyrians, Babylonians, Egyptians, and Greeks. During this time, magnificent temples were constructed, often employing stone as beams. A concise yet insightful history of arches is found, for instance, in Levy \cite{levy}.   Accordingly, the Greeks marked a significant advancement by incorporating {\it curved stones} into more elaborate structures, often for decorative purposes. However, it was the Romans who truly embraced, developed, and mastered the technique of using structural arches to construct immensely intricate monuments such as bridges and aqueducts. In essence, it becomes evident that our ancestors, as they aimed to construct increasingly buildings, needed a more profound understanding of mechanics of the materials  and the
interplay between tension (stress) and deformation (strain).


 Before the fundamental principles  of structural mechanics, constructions probably evolved through a trial-and-error process. With the scientific revolution  around the 16th century, specially in what concerns the infinitesimal differential calculus, the understanding of several laws of nature 
 arose in terms of PDEs systems connected with beam vibrations. 
 More specifically in what concerns  curved bodies like beams, rings, archers, among others,  which are renowned not only for their aesthetic structures but also for their   possible usability, efficiency,  resistance, and  strength, depending on the material composition and  circumstances of applicability, 
  the development of   mathematical models  reflecting structural curved bodies (and their movements)  are increasingly required   in order to make the mechanical process of  tension and deformation more reliable as possible.

In line with the aforementioned context, vibrations of    curved beams and then efforts for {stabilization}  of them over time are subjects that take us back for over at least two centuries ago.
  Among the modelings for curved beams we have nowadays,  it was Bresse (1859) yet in the mid of 19th century that  rigorously derived a reliable 
  set of PDEs for  extensible curved beams under shear and axial effects in the    celebrated work {\it ``Cours de M\'echanique Appliqu\'ee''}   \cite{bresse}. A more complete history about Bresse's discoveries can be seen in the older work by Timoshenko  \cite{timoshenko-history}
 and more recently by Challamel and Elishakoff \cite{challamel-elishakoff}, where the importance of Bresse's contributions are highlighted. We   stress that Timoshenko  (cf. \cite[p. 151]{timoshenko-history}) mentioned the following 
  {\it ``Bresse  was the first
  to take rotatory inertia of the elements of the bar into
  consideration'',} and many other historical data  on Bresse's developments are provided in  \cite[Sect 1]{challamel-elishakoff}. The reliability of Bresse's theory also relies on the fact that other important contributions on circular arches and curved beams take the Bresse development into account, see for instance the remarkable monographs by Timoshenko \cite{timoshenko-part1,timoshenko-part2},  Ciarlet \cite{ciarlet},  Lagnese and Lions \cite{lagnese-lions},   and 
  Lagnese,  Leugering and Schmidt \cite{lag-leug-schm1,lag-leug-schm2}. For this reason, we adopt the Bresse system as the prototype for governing equations in the forthcoming developments of Sections \ref{sec-modeling} and \ref{sec-examplesofsystems}.
  
%

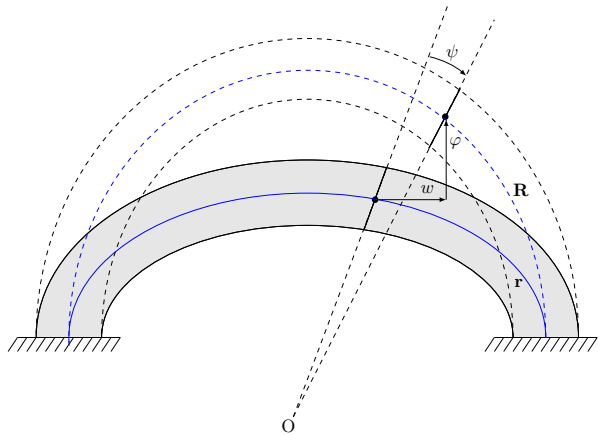
\begin{figure}[htb]
   	\centering
   	\resizebox{8cm}{!}{
   		\centering
   		\begin{tikzpicture}
   			
  \draw[thin,  fill=black!10!white] 
(5.6,0) -- (4.2,0) arc (0:180:4.4cm and 2.4cm) -- (-6.0,0) arc (0:180:-5.8cm and 3.8cm);

  \draw[thin, dashed] (-0.5,-1.7) -- (2.8,7.1);
  \draw[thin, dashed] (-0.5,-1.7) -- (3.8,6.8);
  \node  at (-0.6,-1.9) {O};
    
   \draw[black] (2.4,4.05) -- (3.07,5.33); 
   \draw[black] (2.4,4.05) -- (3.07,5.33);
   \draw[black] (2.4,4.05) -- (3.07,5.33);

   \draw[black] (1.02,2.33) -- (1.5,3.65);
   \draw[black] (1.02,2.33) -- (1.5,3.65);
   \draw[black] (1.02,2.33) -- (1.5,3.65);
   			
    \draw[thin, -latex] (1.26,2.95) -- (2.77,2.95);			
    \fill[black] (1.25,2.95) circle (2pt);
   	\node  at (2.37,3.17) {$w$};
   
   	\draw[thin, -latex] (2.77,2.95) -- (2.77,4.67);
   	\node at (2.97,4.17) {$\varphi$};		
   			
   	\draw[black, latex-] (3.2,5.6) arc (47.8:77.0:1.75cm);
   	\fill[black] (2.75,4.725) circle (2pt);
   	\node  at (2.9, 6.1) {$\psi$};

   			\draw[thin] (-6.4,0) -- (-4.2,0);
   			\draw[thin] (-6.4,0) -- (-6.6,-0.3);
   			\draw[thin] (-6.2,0) -- (-6.4,-0.3);
   			\draw[thin] (-6.0,0) -- (-6.2,-0.3);
   			\draw[thin] (-5.8,0) -- (-6.0,-0.3);
   			\draw[thin] (-5.6,0) -- (-5.8,-0.3);
   			\draw[thin] (-5.4,0) -- (-5.6,-0.3);
   			\draw[thin] (-5.2,0) -- (-5.4,-0.3);
   			\draw[thin] (-5.0,0) -- (-5.2,-0.3);
   			\draw[thin] (-4.8,0) -- (-5.0,-0.3);
   			\draw[thin] (-4.6,0) -- (-4.8,-0.3);
   			\draw[thin] (-4.4,0) -- (-4.6,-0.3);
   			\draw[thin] (-4.2,0) -- (-4.4,-0.3);
   			\draw[thin] (3.8,0) -- (6.0,0);
   			\draw[thin] (3.8,0) -- (3.6,-0.3);
   			\draw[thin] (4.0,0) -- (3.8,-0.3);
   			\draw[thin] (4.2,0) -- (4.0,-0.3);
   			\draw[thin] (4.4,0) -- (4.2,-0.3);
   			\draw[thin] (4.6,0) -- (4.4,-0.3);
   			\draw[thin] (4.8,0) -- (4.6,-0.3);
   			\draw[thin] (5.0,0) -- (4.8,-0.3);
   			\draw[thin] (5.2,0) -- (5.0,-0.3);
   			\draw[thin] (5.4,0) -- (5.2,-0.3);
   			\draw[thin] (5.6,0) -- (5.4,-0.3);
   			\draw[thin] (5.8,0) -- (5.6,-0.3);
   			\draw[thin] (6.0,0) -- (5.8,-0.3);
   			
   \draw[thin] (5.6,0) arc (0:180:5.8cm and 3.8cm);  
   	\draw[thin] (4.2,0) arc (0:180:4.4cm and 2.4cm);
   			
    \draw[thin, blue] (4.9,0) arc (0:182:5.1cm and 3.09cm);
    \node  at (4.32,1.16) {$\mathbf{r}$};

    \draw[thin, dashed] (5.6,0) arc (0:180:5.8cm and 6.4cm);
    \draw[thin, dashed] (4.2,0) arc (0:180:4.4cm and 5.1cm);
 
    \draw[thin, dashed, blue] (4.9,0) arc (0:182:5.1cm and 5.72cm); 
     \node  at (4.37,3.17) {$\mathbf{R}$};
    
   		\end{tikzpicture}
 	}
   	\caption{\small    Longitudinal section  with  an illustrative particle moving from the reference centerline ${\bf r}$ up to the 
    deformed reference line  $\mathbf{R}$  after the displacements $w$, $\varphi$, and $\psi$.}
    \label{fig2}
   \end{figure}

To be  clearer at this point, we note that the well-recognized Bresse system presented in \cite[Chapt. II]{bresse}  (see also Fig. 32 on p. 123 therein for a geometric representation of the curved element and the displacements)  
  provides a very  famous prototype  in (evolution) PDEs that reflects the vibratory movements in   elastic curved bodies under  suitable forces  intrinsically related to {Shear force ($Q$), axial load ($N$), and Bending Moment ($M$)}. Such forces (vibrations)  produce then beam movements such as vertical displacement  ($\varphi$), horizontal displacement ($w$), and angle of rotation ($\psi$)\footnote{In the original work \cite{bresse}, the notations are: $\varphi:=v$, $w:=u$, and $\psi:=\theta$. But here we keep the notation currently used in the literature.}. The  algebraic expressions for the forces $Q,N,M$
 as well as the role  of the displacements $\varphi,w,\psi$  will be  presented in Section \ref{sec-modeling}  with precise details. Below, for the sake of clarity, we provide in Figure \ref{fig2} the geometric illustration of the displacements $\varphi,w,$ and $\psi$.

 In order to control the aforementioned vibrations
 of the elastic   bodies (curved or not), we must direct our attention  towards the stabilization side, which  means the analysis of the asymptotic 
 stability of the  referred displacements. This is closely associated with the concept of damping mechanisms that need to be considered in the stability over time. Among  the possibilities of introduction of internal damping, the present paper is particularly devoted to the thermal one by giving  special attention to the peculiar couplings in systems of Bresse type. 
 Such concern  in asymptotic stability by means of thermoelasticity also  leads us back to studies from at least the mid-20th century, see Dafermos \cite{dafermos} and references therein.
 
Regarding, more specifically, the stability of thermal problems associated with Bresse systems, we can refer several papers from the past two decades, cf. 
  \cite{bazarra-bochicchio-fernandez-naso,bittencourt-silva,bittencourt-silva-camargo, bouraouri-etal23, pedro-hugo-zamp, pedro-hugo-cpaa,delloro-alone,dell-oro-jde2,djellali-labidi-taallah,luci-jaime-IMA,keddi-apalara-messaoudi,liu-rao-zamp09,santos-alone,santos-dilberto-rivera-ADE}, where different constitutive laws for the heat-flux have been considered. As a matter of fact, it is also important to mention that there are several other ways to evaluate the asymptotic dynamics of Bresse systems through different damping mechanisms such as frictional, viscoelastic, and localized effects. See, for instance, \cite{alabou-jaime-dilberto-2011,charles-etal-2013,matofu-rodrigo-siam,santos-soufyane-dilberto,wehbe-youssef}  and references therein.
  The pioneering work on asymptotic stability for curved  beams of Bresse type is given by Liu and Rao \cite{liu-rao-zamp09} where (for the first time) the stabilization properties have been done in terms of the assumption of equal speeds of
  wave propagation
  for a partially damped thermoleastic Bresse model.   As highlighted by the authors in \cite[Sect. 1]{liu-rao-zamp09}, they considered a linearized problem coming from a more general setting in nonlinear thermoelastic flexible beams   as presented in \cite[Chapt. III]{lag-leug-schm2}. 
  In the subsequent years (and until nowadays as one can see in the previous literature)
 several authors have been dealing with different kinds of thermoelastic Bresse systems by showing a diverse range of stability results such as algebraic and exponential decay rates  usually depending on three factors:  proper boundary conditions;   law for the heat flux of conduction; and a relation among the coefficients, often called by {\it stability number} (which reads like {\it equal wave speeds} assumption in Fourier's case).  
 
 Even though there are several  papers dealing with stability for different classes of thermoelastic Bresse models, the authors do not go beyond the stress-strain constitutive laws in order to set the resulting thermoleastic Bresse system.  Indeed, most of them start from Bresse's governing equations and the standard   relations for the shear, axial, and bending forces ($Q,N,M$) in terms of the displacements ($\varphi,\psi,w$) 
 and then connecting the internal thermal component using a straightforward approach  
so that it is not totally clarified how to establish the heat equation and its related coupling terms, because it is usual to  regard the interconnection of  temperature deviations with the forces $Q,N,M$ as a heuristic choice in most cases. As we shall see, thermoelasticity in curved beams of Bresse type possesses an intricate coupling of thermal and mechanical effects, giving rise to distinct designs of coupling terms.
 
Therefore, in Section \ref{sec-modeling}, our first main goal is to explore the fundamental principles underlying thermoelastic effects in curved beams, mainly focusing  on its application to elastic beams of Bresse-type. In this way, we go back to  theoretical foundations of thermo-elasticity 
to comprehend  the 
interplay between temperature variations and the resulting stress-strain relations within  a class of curved beams,  shedding light on the  thermal couplings associated with the   elastic Bresse  system, whose core  is well-known to be expressed in terms of the three   displacements $\varphi$, $\psi$, and $w$,  and their enhanced relations with the forces
$Q$, $N,$ and $M,$ besides the role of temperature and its concise correspondence with  the Bresse system. 
 Our arguments in Section \ref{sec-modeling} are  solidly grounded in the theory of mathematical physics
 provided by   the classical works of
Bresse \cite{bresse},  Ciarlet \cite{ciarlet}, Dafermos \cite{dafermos},  Lagnese,   Leugering and Schmidt \cite{lag-leug-schm1,lag-leug-schm2},  Lagnese and Lions \cite{lagnese-lions}, Timoshenko \cite{timoshenko-part1,timoshenko-part2}, and Timoshenko and Gere \cite{timoshenko-gere}. 

Then,  in Section \ref{sec-examplesofsystems}, our second main goal is to provide new classes of  thermoelastic systems of Bresse type, which are anchored in the solid foundation of Section \ref{sec-modeling} together with distinguished constitutive thermal laws coming from Fourier \cite{fourier}, Cattaneo \cite{cattaneo}, Gurtin and Pipkin \cite{gurtin-pipkin}, Coleman and Gurtin \cite{Coleman-gurtin}, Green and Naghdi \cite{green-naghdi1}, Tzou \cite{tzou95}, Lord-Shulman \cite{lord-shulman}, and still a new class nominated as Bresse-Type III with memory relaxation. In this scenario, a comprehensive list of prototypes of thermoelastic Bresse beam systems come into play as follows:  Bresse-Fourier (BF) models; Bresse-Cattaneo (BC) models;  Bresse-Gurtin-Pipkin (BGP) models;  Bresse-Coleman-Gurtin (BCG) models; Bresse-Green-Naghdi (BGN) models; Bresse-Tzou (BT) models; Bresse-Lord-Shulman (BLS) models; and Bresse-Type III with memory relaxation.

\section{Mathematical modeling} \label{sec-modeling}

    To summarize our development as follows,   we will consider the steps: set the stress tensors upon  both elastic and thermoelastic strains, which   will   be expressed in terms of   displacements and temperature. Then, we will be able to express the forces acting on the governing  Bresse system, since they are known to be represented by the tensors, and finally introduce the representative thermoelastic Bresse systems.

We shall follow initially the framework as developed in \cite{lag-leug-schm1,lag-leug-schm2,lagnese-lions}. Let us start with a thin 3-D beam  of length $L$ with some given initial (small) curvature, $A:=A(x)$ a cross section at a point $x\in[0,L]$, and $\mathbf{r}_0:[0,L]\to \mathbb{R}^3$ the centerline of the beam at rest. With this notation, 
the initial reference configuration of the beam is assumed to occupy the region    $$
\Omega:=\{ \mathbf{r}_0(x_1)+ x_2\mathbf{e}_2+x_3\mathbf{e}_3 ; \ x_1\in[0,L] \mbox{ and } (x_2,x_3)\in A(x_1) \},
 $$ 
which can be understood as the translation of the reference parametric curve $\mathbf{r}_0(x_1)$ to the position $(x_2,x_3):=x_2\mathbf{e}_2+x_3\mathbf{e}_3$  within
the cross section vertical to the tangent of $\mathbf{r}_0$. Our initial assumptions with respect to   cross section and curvature positions are as follows:
The cross section given by $
A(x_1):=\{ (x_2,x_3); x_1\mathbf{e}_1+ x_2\mathbf{e}_2+x_3\mathbf{e}_3\in\Omega \} 
$  is assumed to be  smoothly with respect to $x_1$,   doubly symmetric with respect to  $(x_2,x_3)$  so that
$$
\iint_{A}x_2d(x_2,x_3)=\iint_{A}x_3d(x_2,x_3)=0,
$$
 and with diameter very small when compared to the beam length (thin beams). Also, only initial curvatures are assumed (untwisted beam and  no torsion, neither warping)
	 	so that the simplified skew-symmetric Frenet-type matrix of curvatures is set as
	 	$$
	 	\left(
	 	\begin{array}{ccc}
	 		0 & -\kappa_2 & -\kappa_3	\\
 		\kappa_2 & 0 & 0 \\
 	\kappa_3 & 0 & 0
  	\end{array}\right).
	 	$$

\noindent{\it 1st. Step -- Stress-Strain Relations.} Let us denote by $\sigma_{ij}$ the entries of the stress tensor $\sigma=(\sigma_{ij})_{1\leq i,j\leq3}$. According to the  theory for very thin beams, we may assume that the matrix of stress tensor $\sigma$
is symmetric with $\sigma_{22},\sigma_{23},\sigma_{33}\approx0$, so that it  possesses only three relevant stresses $\sigma_{1j}, \, j=1,2,3$, cf. \cite[Chapt. III, Sect. 2]{lag-leug-schm2}. Then, within the theory of  elastically  and  thermally beams, these remaining stresses  
are given by the following   stress-strain relations  cf. \cite[Chapt. 1, Sect. 6 ]{lagnese-lions}: 
\begin{equation}\label{ssj1}
	\sigma_{1j}=a_{1j}\left(\varepsilon_{1j}-\varepsilon_{1j}^{T}\right), \ \ j=1,2,3,
\end{equation}
 where  the coefficients of elasticity  $a_{1j}$  (to be specified later) are independent of temperature,  the elastic strains $\varepsilon_{1j}$ will be determined below according to precise laws in linear elasticity, and  the notations $\varepsilon_{1j}^{T}$ stand for the  thermal
 strains  whose postulations will be also  theorized  in accordance with proper laws in thermoelastic beams.  Our next steps are to precise the elastic and thermal strains.

 \smallskip

 \noindent{\it 2nd. Step --  Elastic Strains.} In order to provide the strain tensors of the beam in terms of the displacements to be considered  on $\Omega$, let us first outline some notations as follows.

The reference position $\mathbf{r}(x_1,x_2,x_3):=\mathbf{r}_0(x_1)+ x_2\mathbf{e}_2+x_3\mathbf{e}_3$ becomes to $\mathbf{R}(x_1,x_2,x_3,t)$ after deformation (see Figure \ref{fig2}), $\mathbf{V}=\mathbf{R}-\mathbf{r}$ is the displacement vector, and $\mathbf{g}_i:=\mathbf{r}_{x_i}$ is the tangent vector along the undeformed reference parametric curve $\mathbf{r}$. So, the tangent vector after deformation can be given as $\mathbf{G}_i:=\mathbf{g}_i+\mathbf{V}_i, \, i=1,2,3.$ Accordingly, the  Green-St. Venant strain tensor is expressed 
in terms of the displacement gradient  $\varepsilon:=(\nabla \mathbf{V})^T \cdot \nabla \mathbf{V}$, cf. \cite[Sect 1.8]{ciarlet} or \cite[p. 51]{lag-leug-schm2}, and then
$$
\varepsilon_{i,j}:=\frac{1}{2}\left(\frac{\partial\mathbf{R} }{\partial s_i}\cdot \frac{\partial\mathbf{R} }{\partial s_j}-\frac{\partial\mathbf{r}}{\partial s_i} \cdot \frac{\partial\mathbf{r} }{\partial s_j}\right).
$$

By following the first step, the symmetric  matrix of strain tensor $\varepsilon=(\varepsilon_{ij})_{1\leq i,j\leq3}$ is assumed to have   only three  relevant strains, namely, $\varepsilon_{1j}, \, j=1,2,3,$ being $\varepsilon_{22},\varepsilon_{23},\varepsilon_{33}\approx0$ neglected, which is not merely taken for convenience but also it aligns with the theory of thin beams. Thus, using the above notation, we can compute the strains:
\begin{equation}\label{strains1}
\left\{\begin{aligned}
	&\varepsilon_{11} =\frac{1}{2}\left(\frac{\mathbf{G}_1 \cdot \mathbf{G}_1}{ |\mathbf{g}_1|^2}-1\right), \\
	&\varepsilon_{1j} =\frac{1}{2} \frac{1}{|\mathbf{g}_1|}\left(\mathbf{G}_1 \cdot \mathbf{G}_j-\mathbf{g}_1 \cdot \mathbf{g}_j\right), \ j=1,2.
\end{aligned}\right.
\end{equation}
Now, we are going to compute the right above  scalar products
of the Green-St. Venant strain components in terms of  reference strains $\overline{\varepsilon}_{ij}$ and curvatures $\overline{\kappa}_j$, after deformations. To do so, we employ here the particular case of   infinitesimal rigid displacements, which reflects in linear approximations with respect to  $\overline{\varepsilon}_{ij}$, cf. \cite[Sect 6.3]{ciarlet}. Under this situation, the quadratic approximations  
are too small in the sense that they are neglected ($\overline{\varepsilon}_{ij}\overline{\varepsilon}_{kl}\approx0)$.  Yet the remaining reference strains $\overline{\varepsilon}_{22},\overline{\varepsilon}_{23},\overline{\varepsilon}_{33}\approx0$ are usually omitted within the literature of thin beams, so are they here. Thus, still remembering that we do not have  twist nor torsion,  we write down 
\begin{equation}\label{strains2}
\left\{\begin{array}{l}
	\mathbf{G}_1 \cdot \mathbf{G}_1=1+2 \bar{\varepsilon}_{11}-2 x_2\left(\overline{\kappa}_2-2 \bar{\varepsilon}_{12,x_{1}}\right) 
	-2 x_3\left(\overline{\kappa}_3-2 \bar{\varepsilon}_{13,x_{1}}\right), \smallskip \\
	\mathbf{G}_1 \cdot \mathbf{G}_j=2 \bar{\varepsilon}_{1j}, \ j=2,3, \smallskip\\
	|\mathbf{g}_1|\approx1-x_2\kappa_2-x_3\kappa_3, \, \mathbf{g}_1 \cdot \mathbf{g}_j=0,\ j=2,3.
\end{array}\right.
\end{equation}
The identities in \eqref{strains2} can be seen as linear approximations of the quadratic ones obtained in  \cite[Eq. (1.15)]{lag-leug-schm2}. In addition, by setting the final bending strains as
$$
\tilde{\varepsilon}_i:=\overline{\kappa}_i-\kappa_i - 2\overline{\varepsilon}_{1i,x_{1}}, \ i=2,3,
$$
rewriting \eqref{strains2} properly, we can finally consider the linear approximation of the three main Green-St. Venant strains given in \eqref{strains1} as follows
\begin{equation}\label{strain3}
 	\varepsilon_{11}= \bar{\varepsilon}_{11}-x_2 \tilde{\varepsilon}_2-x_3 \tilde{\varepsilon}_3, \ \
	\varepsilon_{12}=  \bar{\varepsilon}_{12}, \ \ 
	\varepsilon_{13}=  \bar{\varepsilon}_{13}.
\end{equation}
Again, the strain identities in \eqref{strain3} can be obtained as a linear (untwisted) approximations of the quadratic 
 ones (finite rotation theory) as reached in  \cite[Eq. (1.17)]{lag-leug-schm2}.  
 In what follows, our goal is to establish  the reference and bending strains ($\bar{\varepsilon}_{1j}$ and $\tilde{\varepsilon}_i$) in terms of the displacements of the curved beam $\Omega$. Accordingly, we can set $\vartheta_1, \vartheta_2,\vartheta_3$ and $W_1,W_2,W_3$ as the rotations and displacements of the beam over the triad $\mathbf{e}_1,\mathbf{e}_2,\mathbf{e}_3$, respectively. Looking for the specific rotations and displacements aligned with Bresse beam theory, we take: 
  \begin{itemize}
  		\item $\vartheta_2:=\psi=\psi(x_1,t)$ as the rotation angle  and
  	$\vartheta_1,\vartheta_3=0$; 
  	
  		\item $W_1:=w=w(x_1,t)$ as the longitudinal (horizontal) displacement; 
  	
 	\item $W_3:=\varphi=\varphi(x_1,t)$ as the  vertical displacement and $W_2=0$.
\end{itemize}
Under the above notation and regarding again the infinitesimal strains and rotations, still in accordance with a linearization of \cite[Eq. (1.18)]{lag-leug-schm2}, we consider the following linear approximations with respect to the elastic reference strains $\overline{\varepsilon}_{1j}$:
\begin{equation}\label{ref-strain}\begin{aligned}
 & \bar{\varepsilon}_{11} := W_{1,x} -\kappa_2 W_2-\kappa_3 W_3 = w_{x_1} -\kappa_3 \varphi,\\ 
 & \bar{\varepsilon}_{12}:=\frac{1}{2}\big(W_{2,x}-\vartheta_3+\kappa_2W_1\big)=\frac{\kappa_2}{2} w, \\ 
 & \bar{\varepsilon}_{13}:=\frac{1}{2}\big(W_{3,x}+\vartheta_2+\kappa_3W_1\big)=\frac{1}{2}(\varphi_{x_1}+\psi+\kappa_3 w).
  	\end{aligned}
\end{equation}

Additionally, with respect to the elastic bending strains $\tilde{\varepsilon}_i$ we firstly follow the approach of \cite[Chapt. IV]{timoshenko-part1} where the stain of longitudinal fibers are proportional to vertical displacements and inversely proportional to the radius of curvature $R$. In such a condition, 
  for small displacements and small curvatures,   the  total curvature difference can be approximated as directly proportional to the vertical displacement as follows
\begin{equation}\label{bending-strain1} \tilde{\varepsilon}_2=\underbrace{\overline{\kappa}_2-\kappa_2}_{\approx\mu E \frac{\varphi}{R}} - 2\overline{\varepsilon}_{12,x_{1}} \approx \mu E \frac{\varphi}{R}  -\kappa_2 w_{x_1},
\end{equation}
where $E$ is  the modulus of elasticity of the beam material and $\mu=\mu(R)$ is a positive constant of proportionality. Second, to the other bending strain, we recall that the deflection of bars with small initial curvature the strain can be related to the rotation of the cross section, as for example assumed in  \cite[Chapt. II]{timoshenko-part2}. See also  the particular situation of   shearable nonlinear 2-d  beams with curvature  \cite[Chapt. III, Sect. 6]{lag-leug-schm2} (which comes from relations in (1.19) therein). Thus, by following this line within our aforementioned assumptions and notations, we put
\begin{equation}\label{bending-strain2}  \tilde{\varepsilon}_3:=\underbrace{\overline{\kappa}_3-\kappa_3}_{\approx (\varphi_{x_1}+\kappa_3w)_{x_1} } - 2\overline{\varepsilon}_{13,x_{1}}\approx- \psi_{x_1}.
\end{equation}

Finally, in order to simplify the notations specially towards Bresse-Timoshenko's framework, we may take
$$
\kappa_2=\kappa_3=\ell>0, \ \  \mu(R)=\frac{1}{R}, \ \ \mbox{ and } \ \ R=\ell^{-1}, 
$$ 
and from \eqref{strain3}, \eqref{ref-strain}, \eqref{bending-strain1} with normalized elasticity ($E=1$), and \eqref{bending-strain2}, we conclude
\begin{equation}\label{strains-final}
\left\{\begin{aligned}
&	\varepsilon_{11}=    (1+x_2\ell)(w_{x_1} -\ell \varphi)+x_3\psi_{x_1},   \\
&	\varepsilon_{12}=  \frac{\ell}{2} w, \\ 
&	\varepsilon_{13}=  \frac{1}{2}(\varphi_{x_1}+\psi+\ell w).
\end{aligned}
\right.	
\end{equation}

\noindent{\it 3rd. Step --  Thermal Strains.}  Apart from considering the elastic displacements, we are still assuming 
that the curved beam is subject to  an  unknown difference of temperature $\Theta(x_1,x_2,x_3,t)$, which significantly contributes to the beam's deformation/regeneration, and its   deviation  is measured in comparison to  a reference state
of uniform temperature distribution $\Theta_{0}(x_1,x_2,x_3)$ in the rest position. In this context,  by following the  principles in  linear thermoelasticity  as in 
\cite[Chapt. I]{lagnese-lions}, for small changes of temperature ($\left|\Theta/\Theta_{0}\right|<<1$) the thermal strains   $\varepsilon_{1j}^{T}$ of \eqref{ssj1} are directly proportional to $\Theta$ and can be represented by 
\begin{equation}\label{therm-strain}
\varepsilon_{1j}^{T}= \alpha \,  \delta_{1 j} \Theta, \quad j=1,2,3,	
\end{equation}
 where $\alpha>0$  is the coefficient  of  {\it thermal  expansion},  $\delta_{1 j}>0$ is a constant of proportionality which depends upon the specific composition of the material comprising the beam. 
 Furthermore, due to the thinness 
 of the beam,  the temperature distribution can be   decomposed as the following {\it Taylor-type}  expansion (cf. \cite[Chapt. III, Sect. 2]{lag-leug-schm2} or else \cite[Sect. 2]{lag-leug-schm1}): 
 $$
 \Theta(x_1,x_2,x_3,t)=\theta^1(x_1,t)+x_2\theta^2(x_1,t)+x_3\theta^3(x_1,t),
 $$
 where $\theta^1, \theta^2, \theta^3$ stand for temperature's components (functions) representing the deviations over $\Omega$. In accordance with  \cite{lag-leug-schm1,lag-leug-schm2,lagnese-lions}, such a condition is  a conventional assumption in the theory of very thin beams, and consequently, we hereby employ  this {\it linearization}  along the cross section  $ A(x_1)$ in our approach. Thus, going back to \eqref{therm-strain}  we obtain the next representations 
 \begin{equation}\label{therm-strain-final}
 	\varepsilon_{1j}^{T}= \alpha \delta_{1 j} \,    \theta^1+x_2 \alpha \delta_{1 j}\theta^2+x_3 \alpha \delta_{1 j}\theta^3 \quad j=1,2,3.	
 \end{equation}

\noindent{\it 4th. Step --  Axial, Bending, and Shear Tensors.} In light of our assumptions following the principles in thermoelastic arched beam vibrations, we consider the standard formulas (cf. \cite[p. 65]{lag-leug-schm2}) for the forces acting over each cross section $A=A(x_1)$  of the $\Omega$ as
\begin{equation*}
\left\{\begin{aligned}
	F_j&=\iint_A \sigma_{1j} \, d(x_2,x_3), \quad j=1,2,3,\\
	M_1&=\iint_A\left(x_2\sigma_{13}-x_3\sigma_{12} \right) d(x_2,x_3), \\
 M_2&=\iint_A x_3 \sigma_{11}   \, d(x_2,x_3), \ \ M_3=-\iint_A x_2 \sigma_{11}  \, d(x_2,x_3).
\end{aligned}\right.
\end{equation*} 

We also set the notations for
area  and inertial moments   of the  cross section $A$, respectively, 
$$
I_1:=\iint_{A} d(x_2,x_3),  \ \   I_2:= \iint_{A} x_2^2 \, d(x_2,x_3),    \ \  I_3:= \iint_{A} x_3^2 \, d(x_2,x_3).
$$ 
Therefore, we can establish the main forces acting over thermoelastic beams of Bresse-Timoshenko type in view of \eqref{ssj1},  \eqref{strains-final}, and \eqref{therm-strain-final}, as it follows. Our main statements come from \cite[Chapts. II and VIII]{bresse} and \cite[Chapt. II]{timoshenko-part2}
\begin{itemize}
	\item {\it $Q$:= Shear Force.} In the Bresse system, the shear force 
is	determined by the equilibrium equation relating the sum of horizontal and vertical forces to the external loads acting on the arch. It is then  influenced by the deformation of the arch due to external loads and the resulting axial force and bending moment. Thus, from the above formulations, the thermo-elastic shear force over the cross section   can be calculated  by means of the following equation:
\begin{equation}\label{shear-componet-temp}
	Q^{\theta^1}:=\iint_A \sigma_{13} \, d(x_2,x_3)\,= \, 
		a_{13}I_1 \bigg[\frac{1}{2}(\varphi_{x_1}+\psi+\ell w) - \alpha\,\delta_{13} \theta^1\bigg].
\end{equation}		
	
	\item {\it $M$:= Bending Moment.} Within classical Bresse-Timoshenko's accomplishment, the bending moment in an arch is determined by the equilibrium equation of statics, which relates the sum of the moments of all the forces about any point on the arch to the external loads acting on the arch. It is also influenced by the deformation of the arch due to external loads and the resulting axial force, being expressed here as the   (thermal-)second bending moment: 
\begin{equation}\label{bending-componet-temp}
	M^{\theta^3}:=\iint_A x_3 \sigma_{11}   \, d(x_2,x_3)\,= \, a_{11}I_3 
	\big[\psi_{x_1} - \alpha\,\delta_{11}   \theta^3\big].
\end{equation}	
		
	\item {\it $N$:= Axial Force.}	With respect to this axial load, we use the formulation where the   principle of equilibrium requires that the sum of the axial and bending loads must be equal to zero, which is a satisfactory approximation for {\it flat} arched beams  under compression (bending) and axial (thrust) stresses, specially those with hinged at the ends (cf. \cite[p. 94-95]{timoshenko-part2}).
Bringing  this peculiar situation to our context,  we may consider $F_1+M_3=0$, and taking into account our notation for axial force $N=F_1$, we finally arrive (without loss of generality\footnote{Indeed, one can see that ${N}^{\theta^2}=\iint_A \sigma_{11} \, d(x_2,x_3)=a_{11}I_1 
	\big[(w_{x_1} -\ell \varphi) - \alpha\, \delta_{11} \theta^1\big]$ which is nothing more than \eqref{Axial-componet-temp} with simplistic mathematical rescalings as $I_1:=I_2\ell$ and $\theta^{1}:=\theta^2 I_2/I_1$. Hence, \eqref{Axial-componet-temp} preserves the core of the axial tensor in terms of elasticity.}) at:	
\begin{equation}\label{Axial-componet-temp}
	{N}^{\theta^2}:=\iint_A x_2 \sigma_{11}  \, d(x_2,x_3)\,= \, a_{11}I_2 
	\big[\ell(w_{x_1} -\ell \varphi) - \alpha\, \delta_{11} \theta^2\big].
\end{equation}		
\end{itemize}

To complete the formulation, one can denote, for instance,  $a_{11}:=E$ and $a_{13}:=2k'G$, where   we recall that $E$ stands for  Young's modulus of elasticity,    $G$ represents the  shear modulus, and  $k'$ means the shear coefficient of correction. Hence, from \eqref{shear-componet-temp}-\eqref{Axial-componet-temp}, the constitutive thermoelastic laws for axial, bending, and shear forces are 
\begin{equation}\label{forces-ABQ}
\left\{\begin{aligned}
Q^{\theta^1}&=k' GI_1 (\varphi_{x_1}+\psi+\ell w) - 2\alpha k' GI_1 \delta_{13}\theta^1, \smallskip \\
M^{\theta^3}&= EI_3\psi_{x_1} - \alpha EI_3  \delta_{11}\theta^3, \smallskip\\
{N}^{\theta^2}&=	E I_2\ell (w_{x_1} -\ell \varphi) - \alpha EI_2 \delta_{11} \theta^2.
\end{aligned}\right.
\end{equation}

It is worth   mentioning that  
the thermoelastic constitutive laws in \eqref{forces-ABQ} completely agree with the classical elastic corresponding  relations when neglecting  thermal effects. Indeed,  one can clearly see that for vanishing thermal strains $\varepsilon_{1j}^T=0$ (i.e. $\Theta=0$ in \eqref{therm-strain} or $\alpha=0$ in \eqref{therm-strain-final}), then 
\eqref{forces-ABQ} with $\theta^j=0, \, j=1,2,3,$ 
reduces to the standard elastic constitutive relations for the axial ($Q^0:=Q$),  
	bending ($M^0:=M$), and shear ($Q^0:=Q$)  forces, namely,
	\begin{equation}\label{ABQ-elastic}
	\left\{\begin{aligned}
		Q&=k'GI_1 (\varphi_{x_1}+\psi+\ell w), \smallskip \\
		M&= EI_3\psi_{x_1} , \smallskip\\
		N&=	E I_2\ell(w_{x_1} -\ell \varphi).
	\end{aligned}\right.
	\end{equation}

\noindent{\it 5th. Step --  Thermodynamic Systems of Bresse type.} Given a general  displacement vector $\mathbf{u}=(u_{i})$ and temperature deviation $T$,   from thermodynamic laws  
the coupled system of partial differential equations that describes the deformation of a material due to thermal and mechanical loads can be  described by the homogeneous system
\begin{equation}\label{thermal-mechanics}
\left\{
\begin{aligned}
 \rho   \partial_{tt} u_i 	&=\big(C_{i j k l} \varepsilon_{kl}	
	\big)_{x_j}-\big(m_{i j} T\big)_{x_j}  \smallskip \\
	 \rho c_\nu{\partial_t T} &= -q_{i,x_i}  -  T_0\, \partial_t \big( m_{i j}\varepsilon_{i j}\big),
\end{aligned}
\right.	
\end{equation}
where  $\rho$ is the density per unit of reference, $c_\nu$ represents the heat capacity, $\varepsilon_{ij}$ are the related strains,  $q_i$ are the components of the heat
flux vector $\mathfrak{q}=(q_i)$, $T_0$ is the reference temperature distribution, and $\left[C_{i j k l}\right],\left[m_{i j}\right]$
are sufficiently smooth tensor fields attributed to the respective reference configuration. For instance, in \cite[Sect. 3]{dafermos} a rigorous justification of \eqref{thermal-mechanics} is given under Fourier's law for the heat flux (e.g. $q_i=-K_{ij} T_{x_j}$) and wave-like strains $\varepsilon_{ij}:=u_{i,x_j}$. See also the motion-energy equations in \cite[Sect. 1]{lord-shulman} and \cite[Chapt. III, Eq. (2.12)]{lag-leug-schm2}  where Fourier's law for heat flux is assumed for the governing flow of heat when the body has a homogeneous reference temperature $T_0$. 
 
 Accordingly, bringing the right above framework to   deformations $(I_1\varphi,I_3\psi,I_1w)$ and thermal ($\Theta$) variations in the body 
 $\Omega,$ and regarding the general  Bresse  principles
  used for the system of vibratory movements (cf. \cite[Chapt. II]{bresse}), we deduce the following  linear thermoelastic governing system 
 of Bresse type
 \begin{equation}\label{bresse-therm-mechanic}
 \left\{\begin{aligned}
  \rho  \partial_{tt}[I_1\varphi]  	&= \big(Q^{\theta^1} + \ell\widehat{{N}^{\theta^2}}\big)_{x_1},  \smallskip \\
 \rho  \partial_{tt}[ I_{3}\psi ] 	& = \big(M^{\theta^3} - \widehat{Q^{\theta^1}}\big)_{x_1}, \smallskip  \\
  \rho  \partial_{tt}[I_1w] 	&= \big({N}^{\theta^2} - \ell\widehat{Q^{\theta^1}}\big)_{x_1}, \smallskip \\
  \rho c_\nu \partial_t\Theta 	&= -\mathbf{q}_{x_1} -  \Theta_{0} \partial_t\big(\varepsilon_{11}+\varepsilon_{13}\big),
 \end{aligned}
 \right.
\end{equation}
where 
$Q^{\theta^1},M^{\theta^3},{N}^{\theta^2}$ are given in  \eqref{forces-ABQ}, the notation $\widehat{F}(x_1)=\int_0^{x_1}F(x)dx$ stands for the primitive,  the strains are given in \eqref{strains-final}, and according to the {\it linearized} expansion for $\Theta(x_1,x_2,x_3,t)  $, the heat flux vector $\mathbf{q}(x_1,x_2,x_3,t)$ is also taken as
$$
\mathbf{q}(x_1,x_2,x_3,t)= q^1(x_1,t)+x_2q^2(x_1,t)+x_3q^3(x_1,t),
$$ 
 where $q^1, q^2, q^3$ are the referred components of the heat flux.

Let us work a  bit more with the fourth equation in \eqref{bresse-therm-mechanic}.  
 Indeed, from the expressions for $\Theta$ and $\mathbf{q}$,  and the strains in \eqref{strains-final},   the   heat  flux of conduction becomes to
 \begin{equation}\label{heat-general2}
 	\begin{aligned}
 		\rho c_\nu \partial_t\Big\{\theta^{1}+x_2\theta^{2}+x_3\theta^{3}\Big\} &   =  -\Big\{q^1_{x_1} +x_2 q^2_{x_1} +x_3q^3_{x_1} \Big\} \\ & 
 		-   \Theta_{0} \partial_t\bigg\{ \frac{1}{2}(\varphi_{x_1}+\psi+\ell w) +(1+x_2\ell)(w_{x_1} -\ell \varphi)+x_3\psi_{x_1}\bigg\}.
 	\end{aligned}
 \end{equation} 
Now,  taking the average (i.e., integrating \eqref{heat-general2} on $A$), the  inertial moments  (i.e., multiplying \eqref{heat-general2} by $x_2$ and $x_3$, respectively, and 
integrating the resulting expression on $A$), and using the reasonable assumption that 
for small deflections  (our current scenario)  the rate of change of the axial force over time can be negligible  when compared to the rate of change of the shear force (cf. \cite[p. 47]{timoshenko-part2}\footnote{For small deflections,  longitudinal forces may be relatively {\it constant} over time as stated therein. Of course the axial force will no longer be constant, but we are assuming its rate of change is negligible compared to   the shear force one in order to  accommodate simplistically the coupling terms,  even metaphysically.}), we are able to split \eqref{heat-general2} into the 
following set of  1D heat equations
$$
\begin{aligned}
	& \rho c_\nu \partial_t\theta^{1}+q^1_{x_1} =- \frac{\Theta_{0}}{2}  (\varphi_{x_1}+\psi+\ell w)_t=-\frac{\Theta_{0}}{2k'GI_1} Q_t:=-\Theta_0^1 Q_t,   \\
	& \rho c_{\nu} \partial_t{\theta}^{2}+q^2_{x_1} 
	=   \Theta_{0}\ell(w_{x_1} -\ell \varphi)_t=-   \frac{\Theta_{0}}{E I_2} N_t :=-\Theta_0^2 N_t,   \smallskip 
	\\	
	& \rho c_{\nu} \partial_t{\theta}^{3}+q^3_{x_1}=- \Theta_{0} \psi_{x_1 t}=-\frac{\Theta_{0}}{EI_3} M_{t}:=\Theta_0^3 M_t,  
\end{aligned}
$$
where in the latter identities we regard the tensors in \eqref{ABQ-elastic}. From this,  the  mechanical and thermal Bresse-type system   \eqref{bresse-therm-mechanic}   can be finally written down   by means of an  expanded class of partial differential equations as follows.

\noindent {\bf I. Full Thermal Coupling:} Bresse system with Shear-Axial-Bending thermal coupling.   
 \begin{equation}\label{bresse-therm-mechanic-final} 
    \left\{
 	\begin{array}{l}
 		\rho I_1 \partial_{tt}\varphi   =  Q^{\theta^1}_{x_1} + \ell{N}^{\theta^2},  \medskip  
 		\\
 		\rho I_{3} \partial_{tt}  \psi  	 =   M^{\theta^3}_{x_1} -  {Q}^{\theta^1} , \medskip 
 		\\
 		\rho I_1 \partial_{tt} w  =  {N}^{\theta^2}_{x_1} - \ell{Q}^{\theta^1}, \medskip
 		\\
 		\rho c_\nu \partial_t\theta^{1} = -q^1_{x_1} -  \Theta_0^1 Q_t, \medskip
 		\\
 		\rho c_{\nu} \partial_t{\theta}^{2} =-q^2_{x_1} 
 		-   \Theta_0^2 N_t, \medskip
 		\\	
 		\rho c_{\nu} \partial_t{\theta}^{3}	=-q^3_{x_1}- \Theta_0^3 M_{t}. 	 
 	\end{array}
 	\right.
 \end{equation}


We remark that \eqref{bresse-therm-mechanic-final} is fully characterized by the triple 
$[Q^{\theta^1},M^{\theta^3},{N}^{\theta^2}]$ given in \eqref{forces-ABQ}, which physically means acknowledging that temperature changes have an impact on the Shear (shearing deformation),  bending (flexural), and axial (tensile/compressive) forces experienced by the arched beam.
Nonetheless, it is reasonable to think that, depending on the body material (among other beam features/stress/components),  temperature variations do not significantly affect the whole structural behavior of the beam, which in turn refers to less  impact of temperature deviations in one (or more) of the distribution of internal forces within the beam. 
This consideration  leads us to a more accurate analysis of the structural response and behavior of the beam under thermal conditions, and yet suggests that deviations of temperature could be diminished (neglected)  
when analyzing their effects on the shear, bending, or/and axial behaviors along the arched beam.
These cases drive us to (sub-)classes of thermoelastic Bresse systems as   described below,  besides the representative system \eqref{bresse-therm-mechanic-final}. In this way, incorporating the following triples with two thermal couplings $ [Q^0,M^{\theta^3},{N}^{\theta^2}]$, $ [Q^{\theta^1},M^0,{N}^{\theta^2}]$, and $[Q^{\theta^1},M^{\theta^3},N^0]$ instead, we  derive the next three classes from \eqref{bresse-therm-mechanic}, respectively. 

\smallskip
 
\noindent {\bf II. Double  Thermal Coupling:}  Bresse systems with  Bending-Axial, Shear-Axial,  and  Shear-Bending thermal coupling, respectively.   
\begin{equation}\label{bresse-therm-mechanic-final-upper}
  \left\{
	\begin{array}{l}
		\rho I_1 \partial_{tt}\varphi   =  Q_{x_1} + \ell{N}^{\theta^2},  \medskip  
		\\
		\rho I_{3} \partial_{tt}  \psi  	 =   M^{\theta^3}_{x_1} -  {Q} , \medskip 
		\\
		\rho I_1 \partial_{tt} w  =  {N}^{\theta^2}_{x_1} - \ell{Q}, \medskip
		\\
		\rho c_{\nu} \partial_t{\theta}^{2} =-q^2_{x_1} 
		-   \Theta_0^2 N_t, \medskip
		\\	
		\rho c_{\nu} \partial_t{\theta}^{3}	=-q^3_{x_1}- \Theta_0^3 M_{t}.	 
	\end{array}
	\right.  \
	 \left\{
	\begin{array}{l}
		\rho I_1 \partial_{tt}\varphi   =  Q^{\theta^1}_{x_1} + \ell{N}^{\theta^2},  \medskip  
		\\
		\rho I_{3} \partial_{tt}  \psi  	 =   M_{x_1} -  {Q}^{\theta^1} , \medskip 
		\\
		\rho I_1 \partial_{tt} w  =  {N}^{\theta^2}_{x_1} - \ell{Q}^{\theta^1}, \medskip
		\\
		\rho c_\nu \partial_t\theta^{1} = -q^1_{x_1} -  \Theta_0^1 Q_t, \medskip
		\\
		\rho c_{\nu} \partial_t{\theta}^{2} =-q^2_{x_1} 
		-   \Theta_0^2 N_t.	 
	\end{array}
	\right.  \  \left\{
	\begin{array}{l}
		\rho I_1 \partial_{tt}\varphi   =  Q^{\theta^1}_{x_1} + \ell{N},  \medskip  
		\\
		\rho I_{3} \partial_{tt}  \psi  	 =   M^{\theta^3}_{x_1} -  {Q}^{\theta^1} , \medskip 
		\\
		\rho I_1 \partial_{tt} w  =  N_{x_1} - \ell{Q}^{\theta^1}, \medskip
		\\
		\rho c_\nu \partial_t\theta^{1} = -q^1_{x_1} -  \Theta_0^1 Q_t, \medskip
		\\
		\rho c_{\nu} \partial_t{\theta}^{3}	=-q^3_{x_1}- \Theta_0^3 M_{t}. 	 
	\end{array}
	\right.
\end{equation}

Last, but no least, by merging   the  triples with one thermal coupling given by $[Q^{\theta^1},M^0,N^0]$, $[Q^0,M^{\theta^3},N^0]$, and $[Q^0,M^0,{N}^{\theta^2}]$, we obtain the following classes of reduced thermoelastic Bresse systems from \eqref{bresse-therm-mechanic}, respectively.

\smallskip
 
\noindent	{\bf III. Single  Thermal Coupling:} Bresse system with    Shear, Bending, and Axial thermal coupling, respectively.     
\begin{equation}\label{bresse-therm-mechanic-final-lower}
\left\{
\begin{array}{l}
	\rho I_1 \partial_{tt}\varphi   =  Q^{\theta^1}_{x_1} + \ell{N},  \medskip  
	\\
	\rho I_{3} \partial_{tt}  \psi  	 =   M_{x_1} -  {Q}^{\theta^1} , \medskip 
	\\
	\rho I_1 \partial_{tt} w  =  N_{x_1} - \ell{Q}^{\theta^1}, \medskip
	\\
	\rho c_\nu \partial_t\theta^{1} = -q^1_{x_1} -  \Theta_0^1 Q_t.  
\end{array}
\right. \
\left\{
\begin{array}{l}
	\rho I_1 \partial_{tt}\varphi   =  Q_{x_1} + \ell{N},  \medskip  
	\\
	\rho I_{3} \partial_{tt}  \psi  	 =   M^{\theta^3}_{x_1} -  {Q} , \medskip 
	\\
	\rho I_1 \partial_{tt} w  =  N_{x_1} - \ell{Q}, \medskip
	\\	
	\rho c_{\nu} \partial_t{\theta}^{3}	=-q^3_{x_1}- \Theta_0^3 M_{t}.	 
\end{array}
\right. \
 \left\{
\begin{array}{l}
	\rho I_1 \partial_{tt}\varphi   =  Q_{x_1} + \ell{N}^{\theta^2},  \medskip  
	\\
	\rho I_{3} \partial_{tt}  \psi  	 =   M_{x_1} -  {Q} , \medskip 
	\\
	\rho I_1 \partial_{tt} w  =  {N}^{\theta^2}_{x_1} - \ell{Q}, \medskip
	\\
	\rho c_{\nu} \partial_t{\theta}^{2} =-q^2_{x_1} 
	-   \Theta_0^2 N_t.
\end{array}
\right.
\end{equation}

%
%

As we are going to clarify in the next section, the  enhanced models provided in   \eqref{bresse-therm-mechanic-final}--\eqref{bresse-therm-mechanic-final-lower}  will furnish a robust tool 
capable of varying several
 thermoelastic Bresse systems.
 Indeed, their versatility lie in the fact of considering different scenarios for thermal couplings (full or partial thermal couplings) and shall incorporate diverse thermal laws for the heat flux components  $q^i,\,i=1,2,3.$  Consequently, we can derive several novel thermoelastic Bresse models, while still providing a deeper justification for those ones previously considered in the  literature.

 \subsection{A short remark on initial-boundary conditions}
 
 In the above development we concentrate our efforts towards the systems of partial differential equations posed on the domain $(0,L)\times(0,\infty)$, with no mention about initial-boundary conditions, namely, conditions at $x=0,L$ and $t=0$. Indeed, before doing so, we initially must pay attention to which model (with full or partial thermal couplings) will come into play when regarding the systems in \eqref{bresse-therm-mechanic-final}--\eqref{bresse-therm-mechanic-final-lower}, since it implies different numbers of variables displayed in the final problem written down.  Besides, to reach a final ruling evolution model including displacements and temperature, we should introduce the constitutive laws for the heat flux components ($q^i,\,i=1,2,3$) in terms of the temperature ones ($\theta^i,\,i=1,2,3$), which in turn will provide several thermoelastic Bresse systems with different physical features that vary a class of problems as previously remarked.  This is exactly the goal of the next section and it still includes models with  memory terms (nonlocal in time) that require not only initial data  but also the past history as well (data for $t\leq0$).

In the existing literature concerned with thermoleastic systems of Bresse type, there are  plenty of examples of boundary conditions that can be taken into account in the forthcoming models, cf. 
 \cite{bazarra-bochicchio-fernandez-naso,bittencourt-silva,bittencourt-silva-camargo, bouraouri-etal23, pedro-hugo-zamp, pedro-hugo-cpaa,delloro-alone,dell-oro-jde2,djellali-labidi-taallah,luci-jaime-IMA,keddi-apalara-messaoudi,lag-leug-schm1,lag-leug-schm2,liu-rao-zamp09,santos-alone,santos-dilberto-rivera-ADE} and references therein.
 However, in order to  orient the reader somehow to possible boundary conditions,  we highlight for instance the physical description of some of them provided in  existing literature. Possible boundary conditions (BC), among many others,  are:
\begin{itemize}
	\item geometric, insulated, or dynamical BC, see e.g. \cite[Chapt. III, Sects. 4-6]{lag-leug-schm2};

	\item some sort of clamped (fixed), hinged (simply supported), or even Cantilever BC, see e.g. \cite[Sect. 2]{liu-lu-curved-bc};
	
		\item  several cases of mixed Dirichlet-Neumann  BC, see e.g. in \cite[Sect. 4]{bittencourt-silva};

	\item the classical (full) Dirichlet BC such as
	$$
	\varphi\Big|_{x=0,L}=\psi\Big|_{x=0,L}= w\Big|_{x=0,L}= \theta^i\Big|_{x=0,L} = q^i\Big|_{x=0,L}=0.
	$$
\end{itemize}

\section{Prototypes of thermoelastic Bresse beams} \label{sec-examplesofsystems}

Our goal in this section is to consider different thermal laws for the heat flux components  $q^i,\,i=1,2,3,$  and then achieve explicitly several thermoelastic beams of Bresse type.
From now on, we observe that we deal with the only one independent spatial variable $x_1$. Consequently, 
we will drop the  subscript $``1"$ by simply representing $x_1$  as $x$, as well as we will use the pattern subscript notations for time and spatial 
derivatives. Additionally, in the next models we adopt the following mathematical notations for the sake of  clarity
\begin{equation}\label{coeff}
	\begin{array}{lll}
		\rho_1=\rho I_1, & k=  k' GI_1, & k_0=EI_1, \smallskip \\
		\rho_2= \rho I_3, & b=EI_3,   &  I_1=I_2\ell, \smallskip \\
		m_1= 2\alpha k'G I_1\delta_{13}, & \varrho_1=\frac{2}{\Theta_0}\rho c_\nu m_1, & \gamma_1=\frac{2}{\Theta_0}m_1, \smallskip \\ 
		m_2= \alpha EI_2\delta_{11}, & \varrho_2=\frac{R}{\Theta_0}\rho c_\nu m_2, & \gamma_2=\frac{R}{\Theta_0}m_2, \smallskip \\
		m_3= \alpha EI_3\delta_{11}, & \varrho_3=\frac{1}{\Theta_0}\rho c_\nu m_3, & \gamma_3=\frac{1}{\Theta_0}m_3. 
	\end{array}  
\end{equation}
Under these notations in \eqref{coeff},  the constitutive laws  \eqref{forces-ABQ}-\eqref{ABQ-elastic}  can be rewritten as
\begin{equation}\label{forces-ABQ-sec3}
\begin{aligned}
	\mbox{\bf Thermoelastic:}	&\left\{\begin{aligned}
		Q^{\theta^1}&=k (\varphi_{x_1}+\psi+\ell w) - m_1\theta^1, \smallskip \\
		M^{\theta^3}&= b\psi_{x_1} - m_3\theta^3, \smallskip\\
		{N}^{\theta^2}&=	k_0 (w_{x_1} -\ell \varphi) - m_2 \theta^2.
	\end{aligned}\right.  \\
	\mbox{\bf Elastic:}	
	& \left\{\begin{aligned}
		Q&=k (\varphi_{x_1}+\psi+\ell w), \smallskip \\
		M&= b\psi_{x_1} , \smallskip\\
		N&=	k_0(w_{x_1} -\ell \varphi).
	\end{aligned}\right.	
\end{aligned}
\end{equation}
Also, the reference models coming from \eqref{bresse-therm-mechanic-final}--\eqref{bresse-therm-mechanic-final-lower} can be simplified  as follows. 

\noindent {\bf I. Full Thermal Coupling:} 
 \begin{equation}\label{bresse-therm-mechanic-full-sec3} 
\left\{
	\begin{array}{l}
		\rho_1  \varphi_{tt}   =  Q^{\theta^1}_{x} + \ell{N}^{\theta^2},  \medskip  
		\\
		\rho_{2}   \psi_{tt}  	 =   M^{\theta^3}_{x} -  {Q}^{\theta^1} , \medskip 
		\\
		\rho_1 w_{tt}  =  N^{\theta^2}_{x} - \ell{Q}^{\theta^1}, \medskip
		\\
		\varrho_1 \theta^{1}_t = - \gamma_1 q^1_{x} - \frac{m_1}{k} Q_t, \medskip
		\\
		\varrho_2 {\theta}^{2}_t =-\gamma_2 q^2_{x} 
		-   \frac{m_2}{k_0} N_t, \medskip
		\\	
		\varrho_3  {\theta}^{3}_t	=-\gamma_3 q^3_{x}-  \frac{m_3}{b} M_{t}. 	 
	\end{array}
	\right.
\end{equation}
\noindent {\bf II. Double Thermal Coupling:}   
\begin{equation}\label{bresse-therm-mechanic-upper-sec3} 
	\left\{
	\begin{array}{l}
		\rho_1 \varphi_{tt}   =  Q_{x} + \ell{N}^{\theta^2},  \medskip  
		\\
		\rho_{2}   \psi_{tt}  	 =   M^{\theta^3}_{x} -  {Q} , \medskip 
		\\
		\rho_1  w_{tt}  =  N^{\theta^2}_{x} - \ell{Q}, \medskip
		\\
	\varrho_2 {\theta}^{2}_t =-\gamma_2 q^2_{x} 
	-   \frac{m_2}{k_0} N_t, \medskip
	\\	
	\varrho_3  {\theta}^{3}_t	=-\gamma_3 q^3_{x}-  \frac{m_3}{b} M_{t}.	 
	\end{array}
	\right.  \
	\left\{
	\begin{array}{l}
		\rho_1 \varphi_{tt}   =  Q^{\theta^1}_{x} + \ell{N}^{\theta^2},  \medskip  
		\\
		\rho_{2}   \psi_{tt}  	 =   M_{x} -  {Q}^{\theta^1} , \medskip 
		\\
		\rho_1  w_{tt}  =  N^{\theta^2}_{x} - \ell{Q}^{\theta^1}, \medskip
		\\
		\varrho_1 \theta^{1}_t = - \gamma_1 q^1_{x} - \frac{m_1}{k} Q_t, \medskip
	\\
	\varrho_2 {\theta}^{2}_t =-\gamma_2 q^2_{x} 
	-   \frac{m_2}{k_0} N_t.	 
	\end{array}
	\right.  \  \left\{
	\begin{array}{l}
		\rho_1 \varphi_{tt}   =  Q^{\theta^1}_{x} + \ell{N},  \medskip  
		\\
		\rho_{2}   \psi_{tt}  	 =   M^{\theta^3}_{x} -  {Q}^{\theta^1} , \medskip 
		\\
		\rho_1  w_{tt}  =  N_{x} - \ell{Q}^{\theta^1}, \medskip
		\\
	\varrho_1 \theta^{1}_t = - \gamma_1 q^1_{x} - \frac{m_1}{k} Q_t, \medskip
		\\
		\varrho_3  {\theta}^{3}_t	=-\gamma_3 q^3_{x}-  \frac{m_3}{b} M_{t}. 	 
	\end{array}
	\right.
\end{equation} 
\noindent 	{\bf III. Single Thermal Coupling:}   
\begin{equation}\label{bresse-therm-mechanic-lower-sec3}
	\left\{
	\begin{array}{l}
		\rho_1 \varphi_{tt}   =  Q^{\theta^1}_{x} + \ell{N},  \medskip  
	\\
	\rho_{2}   \psi_{tt}  	 =   M_{x} -  {Q}^{\theta^1} , \medskip 
	\\
	\rho_1  w_{tt}  =  N_{x} - \ell{Q}^{\theta^1}, \medskip
	\\
	\varrho_1 \theta^{1}_t = - \gamma_1 q^1_{x} - \frac{m_1}{k} Q_t.  
	\end{array}
	\right. \
	\left\{
	\begin{array}{l}
		\rho_1 \varphi_{tt}   =  Q_{x} + \ell{N},  \medskip  
		\\
		\rho_{2}   \psi_{tt}  	 =   M^{\theta^3}_{x} -  {Q} , \medskip 
		\\
		\rho_1  w_{tt}  =  N_{x} - \ell{Q}, \medskip
		\\	
		\varrho_3  {\theta}^{3}_t	=-\gamma_3 q^3_{x}-  \frac{m_3}{b} M_{t}.	 
	\end{array}
	\right. \
	\left\{
	\begin{array}{l}
		\rho_1 \varphi_{tt}   =  Q_{x} + \ell{N}^{\theta^2},  \medskip  
	\\
	\rho_{2}   \psi_{tt}  	 =   M_{x} -  {Q} , \medskip 
	\\
	\rho_1  w_{tt}  =  N^{\theta^2}_{x} - \ell{Q}, \medskip
		\\
		\varrho_2 {\theta}^{2}_t =-\gamma_2 q^2_{x} 
	-   \frac{m_2}{k_0} N_t.
	\end{array}
	\right.
\end{equation}
 
 We are finally in position to state the thermoelastic Bresse models depending on the heat flux of conduction. 

\subsection{Bresse-Fourier (BF) models}

Let us start with the  classical  Fourier's law of heat conduction with a parabolic character, cf. \cite{fourier} (see also 
   \cite{dafermos,lag-leug-schm1,lag-leug-schm2,lagnese-lions} for more details on physical statements): 
\begin{equation}\label{fourier-law}
	q^i= - \varpi_i \theta^i_x,  \ \ i=1,2,3, 
\end{equation}
where $\varpi_i>0$ is a physical parameter. Thus,  replacing \eqref{forces-ABQ-sec3} and \eqref{fourier-law} with $\varpi_i=1$ for simplicity (w.l.o.g.),   the pattern models   \eqref{bresse-therm-mechanic-full-sec3}--\eqref{bresse-therm-mechanic-lower-sec3}, we obtain the next models posed on $(0,L)\times(0,\infty)$.

\smallskip 

\noindent\textbf{I. BF with full thermal coupling.}  
\begin{equation}\label{BF-general}
\left\{	\begin{array}{l}
		\rho_1\varphi_{tt}  - k(\varphi_{x}+\psi + \ell w)_{x} - \ell k_0(w_x - \ell \varphi) + m_1 \theta^1_x + \ell m_2 \theta^2 \, = \, 0,  \medskip \\
		\rho_2\psi_{tt}  -  b\,\psi_{xx} +k(\varphi_{x}+\psi + \ell w) +m_3\theta^3_x - m_1 \theta^1  \, = \, 0, \medskip \\
		\rho_1 w_{tt} - k_0(w_x - \ell \varphi)_x + \ell k(\varphi_{x}+\psi + \ell w) +m_2\theta^2_x- \ell m_1\theta^1 \,= \,0,  \medskip \\
		\varrho_1 \theta^1_t  - \gamma_1 \theta^1_{xx} + m_1(\varphi_x + \psi + \ell w)_t \, = \, 0,   \medskip\\
		\varrho_2 \theta^2_t  - \gamma_2\theta^2_{xx} + m_2(w_x - \ell \varphi)_t \, = \, 0, \medskip \\ 
		\varrho_3 \theta^3_t  - \gamma_3 \theta^3_{xx} + m_3\psi_{xt} \, = \, 0.
	\end{array}
	\right.	
\end{equation}

\noindent\textbf{II. BF with double thermal coupling.} 
\begin{equation}\label{BF-bending-axial}
\hskip-1.4cm  \left\{	\begin{array}{l}
		\rho_1\varphi_{tt}  - k(\varphi_{x}+\psi + \ell w)_{x} - \ell k_0(w_x - \ell \varphi) + \ell m_2 \theta^2 \, = \, 0,  \medskip \\
		\rho_2\psi_{tt}  -  b\,\psi_{xx} +k(\varphi_{x}+\psi + \ell w) +m_3\theta^3_x\, = \, 0, \medskip \\
		\rho_1 w_{tt} - k_0(w_x - \ell \varphi)_x + \ell k(\varphi_{x}+\psi + \ell w) +m_2\theta^2_x \,= \,0,  \medskip \\
		\varrho_2 \theta^2_t  - \gamma_2\theta^2_{xx} + m_2(w_x - \ell \varphi)_t \, = \, 0, \medskip \\ 
		\varrho_3 \theta^3_t  - \gamma_3 \theta^3_{xx} + m_3\psi_{xt} \, = \, 0.
	\end{array}
	\right.	
\end{equation}
\begin{equation}
\label{BF-shear-axial}
 	\left\{	\begin{array}{l}
		\rho_1\varphi_{tt}  - k(\varphi_{x}+\psi + \ell w)_{x} - \ell k_0(w_x - \ell \varphi) + m_1 \theta^1_x + \ell m_2 \theta^2 \, = \, 0,  \medskip \\
		\rho_2\psi_{tt}  -  b\,\psi_{xx} +k(\varphi_{x}+\psi + \ell w) - m_1 \theta^1  \, = \, 0, \medskip \\
		\rho_1 w_{tt} - k_0(w_x - \ell \varphi)_x + \ell k(\varphi_{x}+\psi + \ell w) +m_2\theta^2_x- \ell m_1\theta^1 \,= \,0,  \medskip \\
		\varrho_1 \theta^1_t  - \gamma_1 \theta^1_{xx} + m_1(\varphi_x + \psi + \ell w)_t \, = \, 0,   \medskip\\
		\varrho_2 \theta^2_t  - \gamma_2\theta^2_{xx} + m_2(w_x - \ell \varphi)_t \, = \, 0.
	\end{array}
	\right.	
\end{equation}
\begin{equation}
\label{BF-shear-bending}
\hskip-1.4cm 	\left\{	\begin{array}{l}
		\rho_1\varphi_{tt}  - k(\varphi_{x}+\psi + \ell w)_{x} - \ell k_0(w_x - \ell \varphi) + m_1 \theta^1_x  \, = \, 0,  \medskip \\
		\rho_2\psi_{tt}  -  b\,\psi_{xx} +k(\varphi_{x}+\psi + \ell w) +m_3\theta^3_x - m_1 \theta^1  \, = \, 0, \medskip \\
		\rho_1 w_{tt} - k_0(w_x - \ell \varphi)_x + \ell k(\varphi_{x}+\psi + \ell w) - \ell m_1\theta^1 \,= \,0,  \medskip \\
		\varrho_1 \theta^1_t  - \gamma_1 \theta^1_{xx} + m_1(\varphi_x + \psi + \ell w)_t \, = \, 0,   \medskip\\ 
		\varrho_3 \theta^3_t  - \gamma_3 \theta^3_{xx} + m_3\psi_{xt} \, = \, 0.
	\end{array}
	\right.	
\end{equation}

While the full damped model \eqref{BF-general} is new, and yet not addressed in the literature, the partially damped models \eqref{BF-bending-axial}-\eqref{BF-shear-bending} are strongly aligned with the existing literature. In fact, 
model \eqref{BF-bending-axial} 
stands for
a linear  shearable and flexible thermoelastic beam vibration, which arises as a particular prototype
of nonlinear    thermoelastic flexible beams developed by \cite{lag-leug-schm1,lag-leug-schm2}. Such a model has been addressed by the pioneering\footnote{In the study of  stability for thermoelastic partially damped Bresse systems.} work \cite{liu-rao-zamp09}  and more recently by
\cite{bittencourt-silva}. The stability of \eqref{BF-shear-bending} has been analyzed   in recent times by \cite{bittencourt-silva-camargo} while the stability for \eqref{BF-shear-axial} will appear in a forthcoming work (cf. \cite{bittencourt-silva,bittencourt-silva-camargo}).

\noindent\textbf{III. BF with single thermal coupling.} 
\begin{equation}
	\label{BF-shear}
 	\left\{	\begin{array}{l}
		\rho_1\varphi_{tt}  - k(\varphi_{x}+\psi + \ell w)_{x} - \ell k_0(w_x - \ell \varphi) + m_1 \theta^1_x  \, = \, 0,  \medskip \\
		\rho_2\psi_{tt}  -  b\,\psi_{xx} +k(\varphi_{x}+\psi + \ell w) - m_1 \theta^1  \, = \, 0, \medskip \\
		\rho_1 w_{tt} - k_0(w_x - \ell \varphi)_x + \ell k(\varphi_{x}+\psi + \ell w) - \ell m_1\theta^1 \,= \,0,  \medskip \\
		\varrho_1 \theta^1_t  - \gamma_1 \theta^1_{xx} + m_1(\varphi_x + \psi + \ell w)_t \, = \, 0.  
	\end{array}
	\right.	
\end{equation}
\begin{equation}
\label{BF-bending}
 \hskip-1.6cm	\left\{	\begin{array}{l}
		\rho_1\varphi_{tt}  - k(\varphi_{x}+\psi + \ell w)_{x} - \ell k_0(w_x - \ell \varphi)  \, = \, 0,  \medskip \\
		\rho_2\psi_{tt}  -  b\,\psi_{xx} +k(\varphi_{x}+\psi + \ell w) +m_3\theta^3_x\, = \, 0, \medskip \\
		\rho_1 w_{tt} - k_0(w_x - \ell \varphi)_x + \ell k(\varphi_{x}+\psi + \ell w)  \,= \,0,  \medskip \\
		\varrho_3 \theta^3_t  - \gamma_3 \theta^3_{xx} + m_3\psi_{xt} \, = \, 0.
	\end{array}
	\right.	
\end{equation}
\begin{equation}
\label{BF-axial}
 	\left\{	\begin{array}{l}
		\rho_1\varphi_{tt}  - k(\varphi_{x}+\psi + \ell w)_{x} - \ell k_0(w_x - \ell \varphi) + \ell m_2 \theta^2 \, = \, 0,  \medskip \\
		\rho_2\psi_{tt}  -  b\,\psi_{xx} +k(\varphi_{x}+\psi + \ell w) \, = \, 0, \medskip \\
		\rho_1 w_{tt} - k_0(w_x - \ell \varphi)_x + \ell k(\varphi_{x}+\psi + \ell w) +m_2\theta^2_x \,= \,0,  \medskip \\
		\varrho_2 \theta^2_t  - \gamma_2\theta^2_{xx} + m_2(w_x - \ell \varphi)_t \, = \, 0.
	\end{array}
	\right.	
\end{equation}

Our modeling to achieve the BF models \eqref{BF-shear}, \eqref{BF-bending}, and \eqref{BF-axial},  gives the precise comprehension for the  single thermoelastic coupling   on the shear, bending, and axial, as considered, respectively, by the authors in 
\cite{santos-dilberto-rivera-ADE}, \cite{luci-jaime-IMA}, and \cite{pedro-hugo-zamp} (the latter with Fourier case).

\subsection{Bresse-Cattaneo (BC) models}

Let us take now the Cattaneo \cite{cattaneo} as the thermal relation  (see also \cite{oncu-moodie} for general Cattaneo's type constitutive laws),  which arose as type of hyperbolic heat conduction theory that allows finite speed of wave propagation. It can be seen as
\begin{equation}\label{cattaneo-law}
	 \tau q^i_t +  q^i=- \varpi_i \theta^i_x, \ \ i=1,2,3, 
\end{equation}
where $\tau>0$  is the thermal relaxation time related to the delay in the response of heat flux to a temperature gradient and $\varpi_i>0$. 
Such a theory can be understood as {\it thermoelasticity with one relaxation time.}

Therefore, replacing \eqref{forces-ABQ-sec3} and 
\eqref{cattaneo-law} in the reference models \eqref{bresse-therm-mechanic-full-sec3}--\eqref{bresse-therm-mechanic-lower-sec3}, and choosing $\varpi_i=\gamma_i$ (w.l.o.g.) for the sake of compatibility,  we get the following models on $(0,L)\times(0,\infty)$.

\smallskip 

\noindent\textbf{I. BC with full thermal coupling.} 
\begin{equation}\label{BC-general}
	\left\{	\begin{array}{l}
		\rho_1\varphi_{tt}  - k(\varphi_{x}+\psi + \ell w)_{x} - \ell k_0(w_x - \ell \varphi) + m_1 \theta^1_x + \ell m_2 \theta^2 \, = \, 0,  \medskip \\
		\rho_2\psi_{tt}  -  b\,\psi_{xx} +k(\varphi_{x}+\psi + \ell w) +m_3\theta^3_x - m_1 \theta^1  \, = \, 0, \medskip \\
		\rho_1 w_{tt} - k_0(w_x - \ell \varphi)_x + \ell k(\varphi_{x}+\psi + \ell w) +m_2\theta^2_x- \ell m_1\theta^1 \,= \,0,  \medskip \\
	\varrho_1 \theta^1_t  + \gamma_1 q^1_{x} + m_1(\varphi_x + \psi + \ell w)_t \, = \, 0,   \medskip\\
	\tau q^1_t +  q^1+  \gamma_1 \theta^1_x=0,
\medskip \\
		\varrho_2 \theta^2_t  + \gamma_2q^2_{x} + m_2(w_x - \ell \varphi)_t \, = \, 0, \medskip \\
			\tau q^2_t +  q^2+  \gamma_2 \theta^2_x=0,
		\medskip \\ 
		\varrho_3 \theta^3_t  + \gamma_3q^3_{x} + m_3\psi_{xt} \, = \, 0, 
		\medskip \\
		\tau q^3_t +  q^3+  \gamma_3 \theta^3_x=0.
	\end{array}
	\right.	
\end{equation}
\noindent\textbf{II. BC with double thermal coupling.}
\begin{equation}\label{BC-bending-axial}
\hskip-1.4cm	\left\{	\begin{array}{l}
		\rho_1\varphi_{tt}  - k(\varphi_{x}+\psi + \ell w)_{x} - \ell k_0(w_x - \ell \varphi)  + \ell m_2 \theta^2 \, = \, 0,  \medskip \\
		\rho_2\psi_{tt}  -  b\,\psi_{xx} +k(\varphi_{x}+\psi + \ell w) +m_3\theta^3_x
		\, = \, 0, \medskip \\
		\rho_1 w_{tt} - k_0(w_x - \ell \varphi)_x + \ell k(\varphi_{x}+\psi + \ell w) +m_2\theta^2_x\,= \,0,  \medskip\\
		\varrho_2 \theta^2_t  + \gamma_2q^2_{x} + m_2(w_x - \ell \varphi)_t \, = \, 0, \medskip \\
		\tau q^2_t +  q^2+  \gamma_2 \theta^2_x=0,
		\medskip \\ 
		\varrho_3 \theta^3_t  + \gamma_3q^3_{x} + m_3\psi_{xt} \, = \, 0, 
		\medskip \\
		\tau q^3_t +  q^3+  \gamma_3 \theta^3_x=0.
	\end{array}
	\right.	
\end{equation} 
\begin{equation}\label{BC-shear-axial}
	\left\{	\begin{array}{l}
		\rho_1\varphi_{tt}  - k(\varphi_{x}+\psi + \ell w)_{x} - \ell k_0(w_x - \ell \varphi) + m_1 \theta^1_x + \ell m_2 \theta^2 \, = \, 0,  \medskip \\
		\rho_2\psi_{tt}  -  b\,\psi_{xx} +k(\varphi_{x}+\psi + \ell w)  - m_1 \theta^1  \, = \, 0, \medskip \\
		\rho_1 w_{tt} - k_0(w_x - \ell \varphi)_x + \ell k(\varphi_{x}+\psi + \ell w) +m_2\theta^2_x- \ell m_1\theta^1 \,= \,0,  \medskip \\
		\varrho_1 \theta^1_t  + \gamma_1 q^1_{x} + m_1(\varphi_x + \psi + \ell w)_t \, = \, 0,   \medskip\\
		\tau q^1_t +  q^1+  \gamma_1 \theta^1_x=0,
		\medskip \\
		\varrho_2 \theta^2_t  + \gamma_2q^2_{x} + m_2(w_x - \ell \varphi)_t \, = \, 0, \medskip \\
		\tau q^2_t +  q^2+  \gamma_2 \theta^2_x=0.
	\end{array}
	\right.	
\end{equation} 
\begin{equation}\label{BC-shear-bending}
\hskip-1.4cm	\left\{	\begin{array}{l}
		\rho_1\varphi_{tt}  - k(\varphi_{x}+\psi + \ell w)_{x} - \ell k_0(w_x - \ell \varphi) + m_1 \theta^1_x  \, = \, 0,  \medskip \\
		\rho_2\psi_{tt}  -  b\,\psi_{xx} +k(\varphi_{x}+\psi + \ell w) +m_3\theta^3_x - m_1 \theta^1  \, = \, 0, \medskip \\
		\rho_1 w_{tt} - k_0(w_x - \ell \varphi)_x + \ell k(\varphi_{x}+\psi + \ell w) - \ell m_1\theta^1 \,= \,0,  \medskip \\
		\varrho_1 \theta^1_t  + \gamma_1 q^1_{x} + m_1(\varphi_x + \psi + \ell w)_t \, = \, 0,   \medskip\\
		\tau q^1_t +  q^1+  \gamma_1 \theta^1_x=0,
		\medskip \\
		\varrho_3 \theta^3_t  + \gamma_3q^3_{x} + m_3\psi_{xt} \, = \, 0, 
		\medskip \\
		\tau q^3_t +  q^3+  \gamma_3 \theta^3_x=0.
	\end{array}
	\right.	
\end{equation}

\noindent\textbf{III. BC  with single thermal coupling.} 
\begin{equation}\label{BC-shear}
	\left\{	\begin{array}{l}
		\rho_1\varphi_{tt}  - k(\varphi_{x}+\psi + \ell w)_{x} - \ell k_0(w_x - \ell \varphi) + m_1 \theta^1_x  \, = \, 0,  \medskip \\
		\rho_2\psi_{tt}  -  b\,\psi_{xx} +k(\varphi_{x}+\psi + \ell w)  - m_1 \theta^1  \, = \, 0, \medskip \\
		\rho_1 w_{tt} - k_0(w_x - \ell \varphi)_x + \ell k(\varphi_{x}+\psi + \ell w) - \ell m_1\theta^1 \,= \,0,  \medskip \\
		\varrho_1 \theta^1_t  + \gamma_1 q^1_{x} + m_1(\varphi_x + \psi + \ell w)_t \, = \, 0,   \medskip\\
		\tau q^1_t +  q^1+  \gamma_1 \theta^1_x=0.
	\end{array}
	\right.	
\end{equation}
\begin{equation}\label{BC-bending}
\hskip-1.55cm	\left\{	\begin{array}{l}
		\rho_1\varphi_{tt}  - k(\varphi_{x}+\psi + \ell w)_{x} - \ell k_0(w_x - \ell \varphi)  \, = \, 0,  \medskip \\
		\rho_2\psi_{tt}  -  b\,\psi_{xx} +k(\varphi_{x}+\psi + \ell w) +m_3\theta^3_x\, = \, 0, \medskip \\
		\rho_1 w_{tt} - k_0(w_x - \ell \varphi)_x + \ell k(\varphi_{x}+\psi + \ell w) \,= \,0,  \medskip \\
		\varrho_3 \theta^3_t  + \gamma_3q^3_{x} + m_3\psi_{xt} \, = \, 0, 
		\medskip \\
		\tau q^3_t +  q^3+  \gamma_3 \theta^3_x=0.
	\end{array}
	\right.	
\end{equation} 
\begin{equation}\label{BC-axial}
	\left\{	\begin{array}{l}
		\rho_1\varphi_{tt}  - k(\varphi_{x}+\psi + \ell w)_{x} - \ell k_0(w_x - \ell \varphi)  + \ell m_2 \theta^2 \, = \, 0,  \medskip \\
		\rho_2\psi_{tt}  -  b\,\psi_{xx} +k(\varphi_{x}+\psi + \ell w) 
		\, = \, 0, \medskip \\
		\rho_1 w_{tt} - k_0(w_x - \ell \varphi)_x + \ell k(\varphi_{x}+\psi + \ell w) +m_2\theta^2_x\,= \,0,  \medskip\\
		\varrho_2 \theta^2_t  + \gamma_2q^2_{x} + m_2(w_x - \ell \varphi)_t \, = \, 0, \medskip \\
		\tau q^2_t +  q^2+  \gamma_2 \theta^2_x=0.
	\end{array}
	\right.	
\end{equation}

Some remarks on the BC models are as follows: the modeling of \eqref{BC-bending-axial} gives the justification for the thermoelastic problem presented by  \cite{pedro-hugo-cpaa} (with $\tau>0$ therein); \eqref{BC-bending} justifies a system appearing in \cite{delloro-alone} and later studied by  \cite{keddi-apalara-messaoudi}, as well as the validation of the  Cauchy problem for the BF system  in \cite{said-hamadouche-II} by excluding frictional damping term in the latter ($\gamma=0$ therein). The  remaining BC models have not been approached until the present writing of these notes.

\subsection{Bresse-Gurtin-Pipkin (BGP) models}

With the aim of extending the variety of Bresse's thermoelastic models, let us employ now the  Gurtin-Pipkin thermal law, see Gurtin and Pipkin
\cite{gurtin-pipkin}, where the heat flux of conduction is given by
\begin{equation}\label{gurtinpipkin-law}
 q^i(\cdot,t)= - \varpi_i \int_0^\infty g_i(s) \theta^i_x(\cdot,t-s) ds, \ \ i=1,2,3, 
\end{equation}
where $\varpi_i>0$, the convolutions kernels $g_i:[0, \infty)\to[0, \infty)$ are integrable functions of finite total
mass whose  asymptotic properties will depend on the required study on existence and stability of solution, and the values of $\theta^i$  for negative times are faced as initial data. It is worth mentioning that the Gurtin-Pipkin law \eqref{gurtinpipkin-law} is
more general somehow than the Fourier \eqref{fourier-law} and Cattaneo \eqref{cattaneo-law} ones, see e.g. \cite{conti-pata-squassina-indiana06,delloro-pata-jde-GP}, once for specific choices for the memory kernels $g_i$ we are able to recover them in a rigorous singular limit procedure.  
Accordingly, replacing \eqref{forces-ABQ-sec3} and  \eqref{gurtinpipkin-law} (with $\varpi_i=1$ w.l.o.g.) in the  models 
\eqref{bresse-therm-mechanic-full-sec3}--\eqref{bresse-therm-mechanic-lower-sec3}, we end up with the next thermoelastic BGP models.

\smallskip

\noindent\textbf{I. BGP with full thermal coupling.} 
\begin{equation}\label{BGP-general}
	\left\{	\begin{array}{l}
		\rho_1\varphi_{tt}  - k(\varphi_{x}+\psi + \ell w)_{x} - \ell k_0(w_x - \ell \varphi) + m_1 \theta^1_x + \ell m_2 \theta^2 \, = \, 0,  \medskip \\
		\rho_2\psi_{tt}  -  b\,\psi_{xx} +k(\varphi_{x}+\psi + \ell w) +m_3\theta^3_x - m_1 \theta^1  \, = \, 0, \medskip \\
		\rho_1 w_{tt} - k_0(w_x - \ell \varphi)_x + \ell k(\varphi_{x}+\psi + \ell w) +m_2\theta^2_x- \ell m_1\theta^1 \,= \,0,  \medskip \\
		\varrho_1 \theta^1_t  - \gamma_1 \displaystyle\int_0^\infty g_1(s) \theta^1_{xx}(t-s) ds + m_1(\varphi_x + \psi + \ell w)_t \, = \, 0,   \medskip\\
		\varrho_2 \theta^2_t  - \gamma_2\displaystyle\int_0^\infty g_2(s) \theta^2_{xx}(t-s) ds  + m_2(w_x - \ell \varphi)_t \, = \, 0, \medskip \\ 
		\varrho_3 \theta^3_t  - \gamma_3 \displaystyle\int_0^\infty g_3(s) \theta^3_{xx}(t-s) ds  + m_3\psi_{xt} \, = \, 0.
	\end{array}
	\right.	
\end{equation}

\noindent\textbf{II. BGP with double thermal coupling.}
\begin{equation}\label{BGP-bending-axial}
\hskip-1.45cm	\left\{	\begin{array}{l}
		\rho_1\varphi_{tt}  - k(\varphi_{x}+\psi + \ell w)_{x} - \ell k_0(w_x - \ell \varphi) + \ell m_2 \theta^2 \, = \, 0,  \medskip \\
		\rho_2\psi_{tt}  -  b\,\psi_{xx} +k(\varphi_{x}+\psi + \ell w) +m_3\theta^3_x \, = \, 0, \medskip \\
		\rho_1 w_{tt} - k_0(w_x - \ell \varphi)_x + \ell k(\varphi_{x}+\psi + \ell w) +m_2\theta^2_x \,= \,0,  \medskip \\
		\varrho_2 \theta^2_t  - \gamma_2\displaystyle\int_0^\infty g_2(s) \theta^2_{xx}(t-s) ds  + m_2(w_x - \ell \varphi)_t \, = \, 0, \medskip \\ 
		\varrho_3 \theta^3_t  - \gamma_3 \displaystyle\int_0^\infty g_3(s) \theta^3_{xx}(t-s) ds  + m_3\psi_{xt} \, = \, 0.
	\end{array}
	\right.	
\end{equation} 
\begin{equation}\label{BGP-shear-axial}
	\left\{	\begin{array}{l}
		\rho_1\varphi_{tt}  - k(\varphi_{x}+\psi + \ell w)_{x} - \ell k_0(w_x - \ell \varphi) + m_1 \theta^1_x + \ell m_2 \theta^2 \, = \, 0,  \medskip \\
		\rho_2\psi_{tt}  -  b\,\psi_{xx} +k(\varphi_{x}+\psi + \ell w) - m_1 \theta^1  \, = \, 0, \medskip \\
		\rho_1 w_{tt} - k_0(w_x - \ell \varphi)_x + \ell k(\varphi_{x}+\psi + \ell w) +m_2\theta^2_x- \ell m_1\theta^1 \,= \,0,  \medskip \\
		\varrho_1 \theta^1_t  - \gamma_1 \displaystyle\int_0^\infty g_1(s) \theta^1_{xx}(t-s) ds + m_1(\varphi_x + \psi + \ell w)_t \, = \, 0,   \medskip\\
		\varrho_2 \theta^2_t  - \gamma_2\displaystyle\int_0^\infty g_2(s) \theta^2_{xx}(t-s) ds  + m_2(w_x - \ell \varphi)_t \, = \, 0.
	\end{array}
	\right.	
\end{equation} 
\begin{equation}\label{BGP-shear-bending}
\hskip-1.0cm	\left\{	\begin{array}{l}
		\rho_1\varphi_{tt}  - k(\varphi_{x}+\psi + \ell w)_{x} - \ell k_0(w_x - \ell \varphi) + m_1 \theta^1_x \, = \, 0,  \medskip \\
		\rho_2\psi_{tt}  -  b\,\psi_{xx} +k(\varphi_{x}+\psi + \ell w) +m_3\theta^3_x - m_1 \theta^1  \, = \, 0, \medskip \\
		\rho_1 w_{tt} - k_0(w_x - \ell \varphi)_x + \ell k(\varphi_{x}+\psi + \ell w) - \ell m_1\theta^1 \,= \,0,  \medskip \\
		\varrho_1 \theta^1_t  - \gamma_1 \displaystyle\int_0^\infty g_1(s) \theta^1_{xx}(t-s) ds + m_1(\varphi_x + \psi + \ell w)_t \, = \, 0,   \medskip\\
		\varrho_3 \theta^3_t  - \gamma_3 \displaystyle\int_0^\infty g_3(s) \theta^3_{xx}(t-s) ds  + m_3\psi_{xt} \, = \, 0.
	\end{array}
	\right.	
\end{equation}

\noindent\textbf{III. BGP  with single thermal coupling.} 
\begin{equation}\label{BGP-shear}
	\left\{	\begin{array}{l}
		\rho_1\varphi_{tt}  - k(\varphi_{x}+\psi + \ell w)_{x} - \ell k_0(w_x - \ell \varphi) + m_1 \theta^1_x \, = \, 0,  \medskip \\
		\rho_2\psi_{tt}  -  b\,\psi_{xx} +k(\varphi_{x}+\psi + \ell w)  - m_1 \theta^1  \, = \, 0, \medskip \\
		\rho_1 w_{tt} - k_0(w_x - \ell \varphi)_x + \ell k(\varphi_{x}+\psi + \ell w) - \ell m_1\theta^1 \,= \,0,  \medskip \\
		\varrho_1 \theta^1_t  - \gamma_1 \displaystyle\int_0^\infty g_1(s) \theta^1_{xx}(t-s) ds + m_1(\varphi_x + \psi + \ell w)_t \, = \, 0.
	\end{array}
	\right.	
\end{equation} 
\begin{equation}\label{BGP-bending}
\hskip-2.1cm	\left\{	\begin{array}{l}
		\rho_1\varphi_{tt}  - k(\varphi_{x}+\psi + \ell w)_{x} - \ell k_0(w_x - \ell \varphi)  \, = \, 0,  \medskip \\
		\rho_2\psi_{tt}  -  b\,\psi_{xx} +k(\varphi_{x}+\psi + \ell w) +m_3\theta^3_x  \, = \, 0, \medskip \\
		\rho_1 w_{tt} - k_0(w_x - \ell \varphi)_x + \ell k(\varphi_{x}+\psi + \ell w)  \,= \,0,  \medskip \\
		\varrho_3 \theta^3_t  - \gamma_3 \displaystyle\int_0^\infty g_3(s) \theta^3_{xx}(t-s) ds  + m_3\psi_{xt} \, = \, 0.
	\end{array}
	\right.	
\end{equation} 
\begin{equation}\label{BGP-axial}
\hskip-0.5cm	\left\{	\begin{array}{l}
		\rho_1\varphi_{tt}  - k(\varphi_{x}+\psi + \ell w)_{x} - \ell k_0(w_x - \ell \varphi) + \ell m_2 \theta^2 \, = \, 0,  \medskip \\
		\rho_2\psi_{tt}  -  b\,\psi_{xx} +k(\varphi_{x}+\psi + \ell w) \, = \, 0, \medskip \\
		\rho_1 w_{tt} - k_0(w_x - \ell \varphi)_x + \ell k(\varphi_{x}+\psi + \ell w) +m_2\theta^2_x \,= \,0,  \medskip \\
		\varrho_2 \theta^2_t  - \gamma_2\displaystyle\int_0^\infty g_2(s) \theta^2_{xx}(t-s) ds  + m_2(w_x - \ell \varphi)_t \, = \, 0.
	\end{array}
	\right.	
\end{equation}

Concerning the BGP models \eqref{BGP-general}-\eqref{BGP-axial}, we only have found the papers \cite{delloro-alone,dell-oro-jde2} up to now, where one can see that our models 
  \eqref{BGP-bending} and \eqref{BGP-bending-axial} support physically  the systems addressed in \cite{delloro-alone} and \cite{dell-oro-jde2}, respectively. The remaining BGP models are firstly developed by this work as far as the authors know and, consequently, their stability results (or even {\it numbers of stability}) are still open, although a pattern of stability results could be conjectured.

%
%
%

\subsection{Bresse-Coleman-Gurtin (BCG) models}

An {\it equipresence}  of Fourier and Gurtin-Pipkin thermal laws, yet generating an interesting thermoelastic model,  can be found in the  Coleman-Gurtin heat theory \cite{Coleman-gurtin}, where the resulting model has a parabolic-hyperbolic character due to the connective law 
\begin{equation}\label{coleman-gurtin-law}
	q^i(\cdot,t)= - (1-\varpi_i) \theta^i_x- \varpi_i \int_0^\infty g_i(s) \theta^i_x(\cdot,t-s) ds, \ \ i=1,2,3, 
\end{equation}
with fixed constant $0<\varpi_i<1$  and memory kernel $g_i$ as the previous section. The limit scenarios $\varpi=0,1$ readily match the (parabolic) Fourier \eqref{fourier-law} and the (hyperbolic) Gurtin-Pipkin \eqref{gurtinpipkin-law} laws, respectively. 

The substitution of \eqref{forces-ABQ-sec3} and \eqref{coleman-gurtin-law} in the   models  \eqref{bresse-therm-mechanic-full-sec3}--\eqref{bresse-therm-mechanic-lower-sec3} leads us to the seven thermoelastic BCG systems that very similar to 
\eqref{BGP-general}-\eqref{BGP-axial}. For the sake of brevity, we only consider the most expanded situation \eqref{bresse-therm-mechanic-full-sec3} since the remaining ones are similar to it. Accordingly, in view of \eqref{coleman-gurtin-law} and the forces \eqref{forces-ABQ-sec3}, we obtain the following BCG with full thermal coupling 
\begin{equation}\label{BCG-general}
	\left\{	\begin{array}{l}
		\rho_1\varphi_{tt}  - k(\varphi_{x}+\psi + \ell w)_{x} - \ell k_0(w_x - \ell \varphi) + m_1 \theta^1_x + \ell m_2 \theta^2 \, = \, 0,  \medskip \\
		\rho_2\psi_{tt}  -  b\,\psi_{xx} +k(\varphi_{x}+\psi + \ell w) +m_3\theta^3_x - m_1 \theta^1  \, = \, 0, \medskip \\
		\rho_1 w_{tt} - k_0(w_x - \ell \varphi)_x + \ell k(\varphi_{x}+\psi + \ell w) +m_2\theta^2_x- \ell m_1\theta^1 \,= \,0,  \medskip \\
		\varrho_1 \theta^1_t  -\gamma_1(1-\varpi_1) \theta^1_{xx} - \gamma_1 \varpi_1\displaystyle\int_0^\infty g_1(s) \theta^1_{xx}(t-s) ds + m_1(\varphi_x + \psi + \ell w)_t \, = \, 0,   \medskip\\
		\varrho_2 \theta^2_t -\gamma_2(1-\varpi_2) \theta^2_{xx}  - \gamma_2\varpi_2\displaystyle\int_0^\infty g_2(s) \theta^2_{xx}(t-s) ds  + m_2(w_x - \ell \varphi)_t \, = \, 0, \medskip \\ 
		\varrho_3 \theta^3_t -\gamma_3(1-\varpi_3) \theta^3_{xx}  - \gamma_3 \varpi_3 \displaystyle\int_0^\infty g_3(s) \theta^3_{xx}(t-s) ds  + m_3\psi_{xt} \, = \, 0.
	\end{array}
	\right.	
\end{equation}

The upper and lower partially couplings can be stated similarly. With respect to \eqref{BCG-general}  (and the other BGC partially models) we did not find works on stability on them precisely. However, in \cite{delloro-alone} and \cite{dell-oro-jde2} there are some stability studies for BGP systems where, via a limit procedure, one can reach the  
 same properties for the BCG models related to lower bending moment and upper axial force-bending moment couplings, respectively.

\subsection{Bresse-Green-Naghdi (BGN) models}

In order to exemplify another 
theory of thermoelasticity which possesses  a significant
feature differing from the the classical one (Fourier law), we consider the Green-Naghdi theory \cite{green-naghdi1}, see also \cite{green-naghdi2,green-naghdi3}, where the characterization of material response for
such thermal phenomena relies on three types of constitutive laws, whose terminology is so-called {\it types I, II, and III thermoelasticity}. While the linearized type I thermoelasticiy is strongly related to the classical (parabolic) heat conduction theory, the linearized variants of both  types II and III thermoelasticity grant finite speed of propagation
of thermal waves as in the hyperbolic case, with peculiar lack of thermal damping in the type II case. Roughly speaking, types I and II can be seen as restricted cases of the type III and for this reason we accommodate our BGN models in the third context (then called {\bf BGN-Type III}). 
Accordingly, the type III constitutive law from Green-Naghdi heat flux theory \cite{green-naghdi1,green-naghdi2,green-naghdi3} can be expressed as  
\begin{equation} \label{GN-law}
	q^i=- \beta_i  p^i_x - \varpi_i\theta^i_{x} \  \ \mbox{ with } \ \  p^i(t)=\int_0^t\theta^i(s)ds+p^i(0), \ i=1,2,3,
\end{equation}
 where $p^i$ are called {\it thermal displacements} and  $\beta_i, \varpi_i>0$ are constants related to the thermal conductivity. Note that the threshold cases $\beta_i=0$
 leads us to the classical Fourier law (type I) whereas  \eqref{GN-law} with $\varpi_i=0$ falls into the type II context.   
 
Replacing \eqref{forces-ABQ-sec3} and \eqref{GN-law} in  the models \eqref{bresse-therm-mechanic-full-sec3}--\eqref{bresse-therm-mechanic-lower-sec3} (with simplified coefficients $\beta_i=\beta, \, \varpi_i=\varpi$  w.l.o.g.), we obtain the following BGN-Type III systems.

\smallskip

\noindent\textbf{I. BGN-Type III with full thermal coupling.} 
\begin{equation}\label{BGN-termic-expanded}
	\left\{	\begin{array}{l}
		\rho_1\varphi_{tt}  - k(\varphi_{x}+\psi + \ell w)_{x} - \ell k_0(w_x - \ell \varphi) + m_1 \theta^1_x + \ell m_2 \theta^2 \, = \, 0,  \medskip \\
		\rho_2\psi_{tt}  -  b\,\psi_{xx} +k(\varphi_{x}+\psi + \ell w) +m_3\theta^3_x - m_1 \theta^1  \, = \, 0, \medskip \\
		\rho_1 w_{tt} - k_0(w_x - \ell \varphi)_x + \ell k(\varphi_{x}+\psi + \ell w) +m_2\theta^2_x- \ell m_1\theta^1 \,= \,0,  \medskip \\
		\varrho_1 p^1_{tt}  -\beta \gamma_1 p^1_{xx}-\varpi \gamma_1 p^1_{txx}+ m_1(\varphi_x + \psi + \ell w)_t \, = \, 0,   \medskip\\
			\varrho_2 p^2_{tt}  -\beta \gamma_2 p^2_{xx}-\varpi \gamma_2 p^2_{txx}+ m_2(w_x - \ell \varphi)_t \, = \, 0, \medskip \\ 
			\varrho_3 p^3_{tt}  -\beta \gamma_3 p^3_{xx}-\varpi \gamma_3 p^3_{txx}+ m_3\psi_{xt} \, = \, 0.
	\end{array}
	\right.	
\end{equation}
The last three equations in \eqref{BGN-termic-expanded} can be also rewritten in terms of the temperature components as follows (even on the coming models, but omitted for the sake of reading):
\begin{equation}\label{BGN-heatequations}
	\left\{	\begin{array}{l}
	\varrho_1 \theta^1_{tt}  -\beta \gamma_1 \theta^1_{xx}-\varpi \gamma_1 \theta^1_{txx}+ m_1(\varphi_x + \psi + \ell w)_{tt} \, = \, 0,   \medskip\\
		\varrho_2 \theta^2_{tt}  -\beta \gamma_2 \theta^2_{xx}-\varpi \gamma_2 \theta^2_{txx}+ m_2(w_x - \ell \varphi)_{tt} \, = \, 0, \medskip \\ 
		\varrho_3 \theta^3_{tt}  -\beta \gamma_3 \theta^3_{xx}-\varpi \gamma_3 \theta^3_{txx}+ m_3\psi_{xtt} \, = \, 0.
	\end{array}
	\right.	
\end{equation}

Therefore, the thermoelastic Bresse-Green-Naghdi system of Type III can be seen as the composition of the equations in \eqref{BGN-termic-expanded}$_{1,2,3}$ and \eqref{BGN-heatequations}. The remaining situations are similar and so the details will be neglected.

\noindent\textbf{II. BGN-Type III with double thermal coupling.}
\begin{equation}\label{BGN-bending-axial}
\hskip-1.4cm	\left\{	\begin{array}{l}
		\rho_1\varphi_{tt}  - k(\varphi_{x}+\psi + \ell w)_{x} - \ell k_0(w_x - \ell \varphi) + \ell m_2 \theta^2 \, = \, 0,  \medskip \\
		\rho_2\psi_{tt}  -  b\,\psi_{xx} +k(\varphi_{x}+\psi + \ell w) +m_3\theta^3_x  \, = \, 0, \medskip \\
		\rho_1 w_{tt} - k_0(w_x - \ell \varphi)_x + \ell k(\varphi_{x}+\psi + \ell w) +m_2\theta^2_x \,= \,0,  \medskip \\
	\varrho_2 \theta^2_{tt}  -\beta \gamma_2 \theta^2_{xx}-\varpi \gamma_2 \theta^2_{txx}+ m_2(w_x - \ell \varphi)_{tt} \, = \, 0, \medskip \\ 
\varrho_3 \theta^3_{tt}  -\beta \gamma_3 \theta^3_{xx}-\varpi \gamma_3 \theta^3_{txx}+ m_3\psi_{xtt} \, = \, 0.
	\end{array}
	\right.	
\end{equation} 
\begin{equation}\label{BGN-shear-axial}
	\left\{	\begin{array}{l}
		\rho_1\varphi_{tt}  - k(\varphi_{x}+\psi + \ell w)_{x} - \ell k_0(w_x - \ell \varphi) + m_1 \theta^1_x + \ell m_2 \theta^2 \, = \, 0,  \medskip \\
		\rho_2\psi_{tt}  -  b\,\psi_{xx} +k(\varphi_{x}+\psi + \ell w) - m_1 \theta^1  \, = \, 0, \medskip \\
		\rho_1 w_{tt} - k_0(w_x - \ell \varphi)_x + \ell k(\varphi_{x}+\psi + \ell w) +m_2\theta^2_x- \ell m_1\theta^1 \,= \,0,  \medskip \\
	\varrho_1 \theta^1_{tt}  -\beta \gamma_1 \theta^1_{xx}-\varpi \gamma_1 \theta^1_{txx}+ m_1(\varphi_x + \psi + \ell w)_{tt} \, = \, 0,   \medskip\\
\varrho_2 \theta^2_{tt}  -\beta \gamma_2 \theta^2_{xx}-\varpi \gamma_2 \theta^2_{txx}+ m_2(w_x - \ell \varphi)_{tt} \, = \, 0.
	\end{array}
	\right.	
\end{equation} 
\begin{equation}\label{BGN-shear-bending}
\hskip-1.4cm	\left\{	\begin{array}{l}
		\rho_1\varphi_{tt}  - k(\varphi_{x}+\psi + \ell w)_{x} - \ell k_0(w_x - \ell \varphi) + m_1 \theta^1_x \, = \, 0,  \medskip \\
		\rho_2\psi_{tt}  -  b\,\psi_{xx} +k(\varphi_{x}+\psi + \ell w) +m_3\theta^3_x - m_1 \theta^1  \, = \, 0, \medskip \\
		\rho_1 w_{tt} - k_0(w_x - \ell \varphi)_x + \ell k(\varphi_{x}+\psi + \ell w) - \ell m_1\theta^1 \,= \,0,  \medskip \\
	\varrho_1 \theta^1_{tt}  -\beta \gamma_1 \theta^1_{xx}-\varpi \gamma_1 \theta^1_{txx}+ m_1(\varphi_x + \psi + \ell w)_{tt} \, = \, 0,   \medskip\\ 
\varrho_3 \theta^3_{tt}  -\beta \gamma_3 \theta^3_{xx}-\varpi \gamma_3 \theta^3_{txx}+ m_3\psi_{xtt} \, = \, 0.
	\end{array}
	\right.	
\end{equation}

\noindent\textbf{III. BGN-Type III  with single thermal coupling.} 
\begin{equation}\label{BGN-shear}
	\left\{	\begin{array}{l}
		\rho_1\varphi_{tt}  - k(\varphi_{x}+\psi + \ell w)_{x} - \ell k_0(w_x - \ell \varphi) + m_1 \theta^1_x \, = \, 0,  \medskip \\
		\rho_2\psi_{tt}  -  b\,\psi_{xx} +k(\varphi_{x}+\psi + \ell w)  - m_1 \theta^1  \, = \, 0, \medskip \\
		\rho_1 w_{tt} - k_0(w_x - \ell \varphi)_x + \ell k(\varphi_{x}+\psi + \ell w) - \ell m_1\theta^1 \,= \,0,  \medskip \\
		\varrho_1 \theta^1_{tt}  -\beta \gamma_1 \theta^1_{xx}-\varpi \gamma_1 \theta^1_{txx}+ m_1(\varphi_x + \psi + \ell w)_{tt} \, = \, 0.
	\end{array}
	\right.	
\end{equation} 
\begin{equation}\label{BGN-bending}
\hskip-1.6cm	\left\{	\begin{array}{l}
		\rho_1\varphi_{tt}  - k(\varphi_{x}+\psi + \ell w)_{x} - \ell k_0(w_x - \ell \varphi)  \, = \, 0,  \medskip \\
		\rho_2\psi_{tt}  -  b\,\psi_{xx} +k(\varphi_{x}+\psi + \ell w) +m_3\theta^3_x  \, = \, 0, \medskip \\
		\rho_1 w_{tt} - k_0(w_x - \ell \varphi)_x + \ell k(\varphi_{x}+\psi + \ell w) \,= \,0,  \medskip \\
	\varrho_3 \theta^3_{tt}  -\beta \gamma_3 \theta^3_{xx}-\varpi \gamma_3 \theta^3_{txx}+ m_3\psi_{xtt} \, = \, 0.
	\end{array}
	\right.	
\end{equation} 
\begin{equation}\label{BGN-axial}
	\left\{	\begin{array}{l}
		\rho_1\varphi_{tt}  - k(\varphi_{x}+\psi + \ell w)_{x} - \ell k_0(w_x - \ell \varphi) + \ell m_2 \theta^2 \, = \, 0,  \medskip \\
		\rho_2\psi_{tt}  -  b\,\psi_{xx} +k(\varphi_{x}+\psi + \ell w)   \, = \, 0, \medskip \\
		\rho_1 w_{tt} - k_0(w_x - \ell \varphi)_x + \ell k(\varphi_{x}+\psi + \ell w) +m_2\theta^2_x \,= \,0,  \medskip \\
	\varrho_2 \theta^2_{tt}  -\beta \gamma_2 \theta^2_{xx}-\varpi \gamma_2 \theta^2_{txx}+ m_2(w_x - \ell \varphi)_{tt} \, = \, 0.
	\end{array}
	\right.	
\end{equation}

While the full damped system \eqref{BGN-termic-expanded} and the upper partially damped systems \eqref{BGN-bending-axial}-\eqref{BGN-shear-bending} have not been studied in the literature so far, the lower thermal couplings generating \eqref{BGN-shear}, \eqref{BGN-bending}, and \eqref{BGN-axial}, give the precise physical justification for the core of the thermoelastic systems appearing in \cite{santos-alone}, \cite{said-hamadouche}, and \cite{bouraouri-etal23}, respectively. We note that in the works \cite{bouraouri-etal23,said-hamadouche}  the authors added mathematically frictional-viscous damping terms in order to evaluate the stability of the related systems therein. We still stress that \eqref{BGN-bending}  justifies precisely a model appearing  in \cite{djellali-labidi-taallah} in its  thermal essence but  the  authors have changed  indiscriminately the Kelvin-Voigt damping term in \eqref{BGN-bending}$_{4}$ to a historyless memory term with the statement that their heat equation {\it ``represents a heat conduction problem with memory in the sense of Green
	and Naghdi, it models the dependence on the history of the temperature gradient''\footnote{See on p. 3 therein.}.} As we are going to see further, such a model considered in \cite{djellali-labidi-taallah}  can be ultimately obtained through a rigorous new heat-flux law invoking  type III with memory relaxation.

%
%
 
\subsection{Bresse-Tzou (BT) models}

To amplify the variety of thermoelastic Bresse systems yielding from the governing equations in the previous section along with different constitutive thermal laws, we consider now the so-called dual-phase-lag
heat conduction theory firstly proposed by Tzou \cite{tzou95}
where two delays come  into play in the constitutive equations. Here, we address the Taylor approximations of   second order for the heat flux and of  first order   for the temperature, namely,  
\begin{equation}\label{Tzou-law}
{q}^i+\tau_q {q}^i_t+\frac{\tau_q^2}{2} {q}^i_{tt}=-\varpi \theta^i_{x}- \varpi  \tau_\theta \theta^i_{tx}, \ \ i=1,2,3,	
\end{equation}
where the time delays $\tau_q$ and $\tau_\theta$ are related to relaxation
time and microstructural interactions, respectively, and $\varpi>0$. In the context of Bresse beam systems,  the dual-phase-lag heat conduction model may be employed to capture more accurate temperature distributions, by considering the finite speed of heat propagation (hyperbolic character) and two different relaxation times. This is why such a theory can be called {\it thermoelasticity with two relaxation times.} We still note that neglecting the second approximation for the heat flux ($\tau^2_q\approx0$) and the first approximation for the temperature $(\tau_\theta\approx0$), then \eqref{Tzou-law} becomes the Cattaneo law \eqref{cattaneo-law}.  In such a case, the coming procedure will cover the Lord-Shulman theory \cite{lord-shulman} and, therefore, we have another class of thermoelastic problems called thereafter by Bresse-Lord-Shulman (BLS) systems. 
In this way, the  dual-phase-lag BT models (here called {\bf BT-DPL}) 
provide, let us say, a more  sophisticated description of the thermal behavior in Bresse beam systems in what concerns ultrafast processes.

Again from \eqref{forces-ABQ-sec3} and relying on \eqref{Tzou-law} (with $\varpi=1$ for simplicity), then the reference models \eqref{bresse-therm-mechanic-full-sec3}--\eqref{bresse-therm-mechanic-lower-sec3}  turn  into the following BT-DPL thermoelastic systems\footnote{The convenient notation $\partial_t^j$ meaning the $j$-th  derivative has been employed.}. 

\smallskip 

\noindent\textbf{I. BT-DPL with full thermal coupling.} 
\begin{equation}\label{BT-termic-expanded}
	\left\{	\begin{array}{l}
		\rho_1\varphi_{tt}  - k(\varphi_{x}+\psi + \ell w)_{x} - \ell k_0(w_x - \ell \varphi) + m_1 \theta^1_x + \ell m_2 \theta^2 \, = \, 0,  \medskip \\
		\rho_2\psi_{tt}  -  b\,\psi_{xx} +k(\varphi_{x}+\psi + \ell w) +m_3\theta^3_x - m_1 \theta^1  \, = \, 0, \medskip \\
		\rho_1 w_{tt} - k_0(w_x - \ell \varphi)_x + \ell k(\varphi_{x}+\psi + \ell w) +m_2\theta^2_x- \ell m_1\theta^1 \,= \,0,  \medskip \\
	\varrho_1\left[\frac{\tau_q^2}{2} \theta^1_{t t t}+\tau_q \theta^1_{t t}+\theta^1_t\right]- \gamma_1\left[\theta^1_{xx}+ \tau_\theta  \theta^1_{txx}\right]+m_1(1+\tau_q\partial_t+\frac{\tau_q^{2}}{2}\partial_{tt}) (\varphi_x + \psi + \ell w)_t=0,   \medskip\\
	\varrho_2\left[\frac{\tau_q^2}{2} \theta^2_{t t t}+\tau_q \theta^2_{t t}+\theta^2_t\right]- \gamma_2\left[\theta^2_{xx}+ \tau_\theta  \theta^2_{txx}\right]+m_2(1+\tau_q\partial_t+\frac{\tau_q^{2}}{2}\partial_{tt})
	(w_x - \ell \varphi)_{t}=0, \medskip \\ 
	\varrho_3\left[\frac{\tau_q^2}{2} \theta^3_{t t t}+\tau_q \theta^3_{t t}+\theta^3_t\right]- \gamma_3\left[\theta^3_{xx}+ \tau_\theta  \theta^3_{txx}\right]+m_3(1+\tau_q\partial_t+\frac{\tau_q^{2}}{2}\partial_{tt})\psi_{xt}=0.
	\end{array}
	\right.	
\end{equation}

\noindent\textbf{II. BT-DPL  with double thermal coupling.} 
\begin{equation}\label{BT-bending-axial}
\hskip-0.7cm	\left\{	\begin{array}{l}
		\rho_1\varphi_{tt}  - k(\varphi_{x}+\psi + \ell w)_{x} - \ell k_0(w_x - \ell \varphi) + \ell m_2 \theta^2 \, = \, 0,  \medskip \\
		\rho_2\psi_{tt}  -  b\,\psi_{xx} +k(\varphi_{x}+\psi + \ell w) +m_3\theta^3_x   \, = \, 0, \medskip \\
		\rho_1 w_{tt} - k_0(w_x - \ell \varphi)_x + \ell k(\varphi_{x}+\psi + \ell w) +m_2\theta^2_x\,= \,0,  \medskip \\
		\varrho_2\left[\frac{\tau_q^2}{2} \theta^2_{t t t}+\tau_q \theta^2_{t t}+\theta^2_t\right]- \gamma_2\left[\theta^2_{xx}+ \tau_\theta  \theta^2_{txx}\right]+m_2(1+\tau_q\partial_t+\frac{\tau_q^{2}}{2}\partial_{tt})
		(w_x - \ell \varphi)_{t}=0, \medskip \\ 
		\varrho_3\left[\frac{\tau_q^2}{2} \theta^3_{t t t}+\tau_q \theta^3_{t t}+\theta^3_t\right]- \gamma_3\left[\theta^3_{xx}+ \tau_\theta  \theta^3_{txx}\right]+m_3(1+\tau_q\partial_t+\frac{\tau_q^{2}}{2}\partial_{tt})\psi_{xt}=0.
	\end{array}
	\right.	
\end{equation} 
\begin{equation}\label{BT-shear-axial}
	\left\{	\begin{array}{l}
		\rho_1\varphi_{tt}  - k(\varphi_{x}+\psi + \ell w)_{x} - \ell k_0(w_x - \ell \varphi) + m_1 \theta^1_x + \ell m_2 \theta^2 \, = \, 0,  \medskip \\
		\rho_2\psi_{tt}  -  b\,\psi_{xx} +k(\varphi_{x}+\psi + \ell w) - m_1 \theta^1  \, = \, 0, \medskip \\
		\rho_1 w_{tt} - k_0(w_x - \ell \varphi)_x + \ell k(\varphi_{x}+\psi + \ell w) +m_2\theta^2_x- \ell m_1\theta^1 \,= \,0,  \medskip \\
		\varrho_1\left[\frac{\tau_q^2}{2} \theta^1_{t t t}+\tau_q \theta^1_{t t}+\theta^1_t\right]- \gamma_1\left[\theta^1_{xx}+ \tau_\theta  \theta^1_{txx}\right]+m_1(1+\tau_q\partial_t+\frac{\tau_q^{2}}{2}\partial_{tt}) (\varphi_x + \psi + \ell w)_t=0,   \medskip\\
		\varrho_2\left[\frac{\tau_q^2}{2} \theta^2_{t t t}+\tau_q \theta^2_{t t}+\theta^2_t\right]- \gamma_2\left[\theta^2_{xx}+ \tau_\theta  \theta^2_{txx}\right]+m_2(1+\tau_q\partial_t+\frac{\tau_q^{2}}{2}\partial_{tt})
		(w_x - \ell \varphi)_{t}=0.
	\end{array}
	\right.	
\end{equation} 
\begin{equation}\label{BT-shear-bending}
	\left\{	\begin{array}{l}
		\rho_1\varphi_{tt}  - k(\varphi_{x}+\psi + \ell w)_{x} - \ell k_0(w_x - \ell \varphi) + m_1 \theta^1_x \, = \, 0,  \medskip \\
		\rho_2\psi_{tt}  -  b\,\psi_{xx} +k(\varphi_{x}+\psi + \ell w) +m_3\theta^3_x - m_1 \theta^1  \, = \, 0, \medskip \\
		\rho_1 w_{tt} - k_0(w_x - \ell \varphi)_x + \ell k(\varphi_{x}+\psi + \ell w) - \ell m_1\theta^1 \,= \,0,  \medskip \\
		\varrho_1\left[\frac{\tau_q^2}{2} \theta^1_{t t t}+\tau_q \theta^1_{t t}+\theta^1_t\right]- \gamma_1\left[\theta^1_{xx}+ \tau_\theta  \theta^1_{txx}\right]+m_1(1+\tau_q\partial_t+\frac{\tau_q^{2}}{2}\partial_{tt}) (\varphi_x + \psi + \ell w)_t=0,   \medskip\\ 
		\varrho_3\left[\frac{\tau_q^2}{2} \theta^3_{t t t}+\tau_q \theta^3_{t t}+\theta^3_t\right]- \gamma_3\left[\theta^3_{xx}+ \tau_\theta  \theta^3_{txx}\right]+m_3(1+\tau_q\partial_t+\frac{\tau_q^{2}}{2}\partial_{tt})\psi_{xt}=0.
	\end{array}
	\right.	
\end{equation}

\noindent\textbf{III. BT-DPL   with single thermal coupling.} 
\begin{equation}\label{BT-shear}
	\left\{	\begin{array}{l}
		\rho_1\varphi_{tt}  - k(\varphi_{x}+\psi + \ell w)_{x} - \ell k_0(w_x - \ell \varphi) + m_1 \theta^1_x \, = \, 0,  \medskip \\
		\rho_2\psi_{tt}  -  b\,\psi_{xx} +k(\varphi_{x}+\psi + \ell w) - m_1 \theta^1  \, = \, 0, \medskip \\
		\rho_1 w_{tt} - k_0(w_x - \ell \varphi)_x + \ell k(\varphi_{x}+\psi + \ell w) - \ell m_1\theta^1 \,= \,0,  \medskip \\
		\varrho_1\left[\frac{\tau_q^2}{2} \theta^1_{t t t}+\tau_q \theta^1_{t t}+\theta^1_t\right]- \gamma_1\left[\theta^1_{xx}+ \tau_\theta  \theta^1_{txx}\right]+m_1  (1+\tau_q\partial_t+\frac{\tau_q^{2}}{2}\partial_{tt}) (\varphi_x + \psi + \ell w)_t =0.
	\end{array}
	\right.	
\end{equation} 
\begin{equation}\label{BT-bending}
\hskip-3.1cm	\left\{	\begin{array}{l}
		\rho_1\varphi_{tt}  - k(\varphi_{x}+\psi + \ell w)_{x} - \ell k_0(w_x - \ell \varphi) \, = \, 0,  \medskip \\
		\rho_2\psi_{tt}  -  b\,\psi_{xx} +k(\varphi_{x}+\psi + \ell w) +m_3\theta^3_x \, = \, 0, \medskip \\
		\rho_1 w_{tt} - k_0(w_x - \ell \varphi)_x + \ell k(\varphi_{x}+\psi + \ell w)  \,= \,0,  \medskip \\
		\varrho_3\left[\frac{\tau_q^2}{2} \theta^3_{t t t}+\tau_q \theta^3_{t t}+\theta^3_t\right]- \gamma_3\left[\theta^3_{xx}+ \tau_\theta  \theta^3_{txx}\right]+m_3 (1+\tau_q\partial_t+\frac{\tau_q^{2}}{2}\partial_{tt})\psi_{xt} =0.
	\end{array}
	\right.	
\end{equation}  
\begin{equation}\label{BT-axial}
\hskip-0.8cm	\left\{	\begin{array}{l}
		\rho_1\varphi_{tt}  - k(\varphi_{x}+\psi + \ell w)_{x} - \ell k_0(w_x - \ell \varphi) + \ell m_2 \theta^2 \, = \, 0,  \medskip \\
		\rho_2\psi_{tt}  -  b\,\psi_{xx} +k(\varphi_{x}+\psi + \ell w)   \, = \, 0, \medskip \\
		\rho_1 w_{tt} - k_0(w_x - \ell \varphi)_x + \ell k(\varphi_{x}+\psi + \ell w) +m_2\theta^2_x\,= \,0,  \medskip \\
		\varrho_2\left[\frac{\tau_q^2}{2} \theta^2_{t t t}+\tau_q \theta^2_{t t}+\theta^2_t\right]- \gamma_2\left[\theta^2_{xx}+ \tau_\theta  \theta^2_{txx}\right]+m_2(1+\tau_q\partial_t+\frac{\tau_q^{2}}{2}\partial_{tt})
		(w_x - \ell \varphi)_{t}=0.
	\end{array}
	\right.	
\end{equation}

The set of  BT-DPL models \eqref{BT-termic-expanded}-\eqref{BT-axial}  provide a new  class of thermoelastic Bresse systems with dual-phase-lag
heat conduction, and  yet not considered in the literature  except for system \eqref{BT-bending}. As a matter of fact, starting from \eqref{BT-bending} and introducing the following notation  $\hat{f}:=f+\tau_q f_t+ \frac{\tau_q^2}{2} f_{tt}$, we are able to reach precisely (and justify physically) the BT-DPL system considered very recently by \cite{bazarra-bochicchio-fernandez-naso}, where numerical studies are accounted for stability. In addition, we remark that  theoretical  results to achieve stability properties for the aforementioned BT-DPL models, a new kind of {\it stability number} involving the coefficients of the system will certainly  be required, besides the condition imposed on the delay parameters  $\tau_q$ and $\tau_\theta$ (cf. \cite{bazarra-bochicchio-fernandez-naso}).

%
%
%
%
%

\subsection{Bresse-Lord-Shulman (BLS) models}

As remarked previously, 
 the 
dual-phase-lag
heat conduction  with only first-order
approximation for the heat flux still covers the hyperbolic Lord-Shulmann (cf. \cite{lord-shulman}) model. In other words, it can be done by 
taking \eqref{Tzou-law} with vanishing  second-order heat flux ($\tau^2_q\approx0$) and first-order temperature gradient ($\tau_\theta\approx0$) phase legs. Therefore, proceeding verbatim as done for the BT-DPL models \eqref{BT-termic-expanded}-\eqref{BT-axial}, we can state a class of thermoelastic BLS systems. For the sake of shortness, we only consider the case related to   \eqref{BT-termic-expanded} and write down the BLS system with 
 full thermal coupling only (yet denoting $\tau_q:=\tau$). 
\begin{equation}\label{BLS-termic-expanded}
	\left\{	\begin{array}{l}
		\rho_1\varphi_{tt}  - k(\varphi_{x}+\psi + \ell w)_{x} - \ell k_0(w_x - \ell \varphi) + m_1 \theta^1_x + \ell m_2 \theta^2 \, = \, 0,  \medskip \\
		\rho_2\psi_{tt}  -  b\,\psi_{xx} +k(\varphi_{x}+\psi + \ell w) +m_3\theta^3_x - m_1 \theta^1  \, = \, 0, \medskip \\
		\rho_1 w_{tt} - k_0(w_x - \ell \varphi)_x + \ell k(\varphi_{x}+\psi + \ell w) +m_2\theta^2_x- \ell m_1\theta^1 \,= \,0,  \medskip \\
		\varrho_1(\tau  \theta^1_{t t}+\theta^1_t)- \gamma_1 \theta^1_{xx} +m_1 (1+\tau\partial_t)(\varphi_x + \psi + \ell w)_t  =0,   \medskip\\
		\varrho_2(\tau  \theta^2_{t t}+\theta^2_t)- \gamma_2 \theta^2_{xx}+m_2(1+\tau\partial_t)(w_x - \ell \varphi)_{t}  =0, \medskip \\ 
		\varrho_3(\tau  \theta^3_{t t}+\theta^3_t)- \gamma_3 \theta^3_{xx} +m_3 (1+\tau\partial_t)\psi_{xt}  =0.
	\end{array}
	\right.	
\end{equation}
 
The remaining cases of partially damped BLS systems can be done similarly and there is no knowledge up to now about {\it stability number} for the analysis of their stability properties.

\subsection{Bresse-Type III with memory relaxation}

To finish this work, we are going to propose a new thermoelastic constitutive law based on the harmonization of the  Green-Naghdi  \cite{green-naghdi1} and Gurtin-Pipkin
\cite{gurtin-pipkin} theories. Indeed, revisiting  the type III heat flow section in \cite[Sect. 8.3]{green-naghdi1}
where the analysis is analogously made upon the flow of ``viscoelastic" bodies with heat-flux given in terms of the thermal displacement and temperature gradient, but yet regarding the linearized theory for infinitesimal temperature gradient of \cite[Sect. 7]{gurtin-pipkin}
where the thermal functional related to  heat-flux  is derived from the Riesz Theorem as a weighted-integral operator, then the following
response heat flux vectors shall be accommodated as the melding of the thermal displacement and viscoelastic relaxation (in the context of memory) of the temperature gradient. Accordingly, we deduce the constitutive equation for $i=1,2,3:$ 
\begin{equation} \label{MM-law}
	q^i=- \beta_i  p^i_x - \varpi_i \int_0^\infty \mu_i(s) \theta^i_x(t-s) ds, \  \ \mbox{ with } \ \  p^i_t(t)=\theta^i(t), \ t\in \mathbb{R}.
\end{equation}
 As previously, $p^i$ are the {thermal displacements},  $\beta_i, \varpi_i>0$, and $\mu_i$ are differentiable functions known as {\it heat-flux relaxation  kernels}, with $\mu_i(0)$ the instantaneous conductivity, and $\mu_i(\infty)=0.$
 Equation \eqref{MM-law} can also be seen as a memory relaxation of the Green-Naghdi law \eqref{GN-law},  reason 
 why we designate it as 
{\it type III constitutive law with memory relaxation} for the heat flux. It is worth pointing out that it induces a fully hyperbolic character of heat transfer, similar to both Gurtin-Pipkin and Green-Naghdi laws. 

Finally, substituting \eqref{forces-ABQ-sec3} and the constitutive law \eqref{MM-law} in the models \eqref{bresse-therm-mechanic-full-sec3}--\eqref{bresse-therm-mechanic-lower-sec3}
(again with $\beta_i=\beta, \, \varpi_i=\varpi$  w.l.o.g.), we get the following  Bresse systems with Memory Relaxation of Type III,    hereby named as {\bf BMR-TypeIII}.

\smallskip

\noindent\textbf{I. BMR-TypeIII  with full thermal coupling.} 
\begin{equation}\label{MM-termic-full}
	\left\{	\begin{array}{l}
		\rho_1\varphi_{tt}  - k(\varphi_{x}+\psi + \ell w)_{x} - \ell k_0(w_x - \ell \varphi) + m_1 \theta^1_x + \ell m_2 \theta^2 \, = \, 0,  \medskip \\
		\rho_2\psi_{tt}  -  b\,\psi_{xx} +k(\varphi_{x}+\psi + \ell w) +m_3\theta^3_x - m_1 \theta^1  \, = \, 0, \medskip \\
		\rho_1 w_{tt} - k_0(w_x - \ell \varphi)_x + \ell k(\varphi_{x}+\psi + \ell w) +m_2\theta^2_x- \ell m_1\theta^1 \,= \,0,  \medskip \\
		\varrho_1 p^1_{tt}  -\beta \gamma_1 p^1_{xx}
		-\varpi \gamma_1\displaystyle\int_0^\infty \mu_1(s) p^1_{txx}(t-s) ds + m_1(\varphi_x + \psi + \ell w)_t \, = \, 0,   \medskip\\
		\varrho_2 p^2_{tt}  -\beta \gamma_2 p^2_{xx}-\varpi \gamma_2 \displaystyle\int_0^\infty \mu_2(s) p^2_{txx}(t-s) ds+ m_2(w_x - \ell \varphi)_t \, = \, 0, \medskip \\ 
		\varrho_3 p^3_{tt}  -\beta \gamma_3 p^3_{xx}
			-\varpi \gamma_3 \displaystyle\int_0^\infty \mu_3(s) p^3_{txx}(t-s) ds+ m_3\psi_{xt} \, = \, 0.
	\end{array}
	\right.	
\end{equation}

This approach   allows us to reach mathematically  the last three equations in the  previous BGN-Type III  model
\eqref{BGN-termic-expanded} by
simply proceeding as, e.g.,  in \cite{CDGP-I,CDGP-II}. Indeed, by taking the integrated kernels $\mu_i(s)$ converging (in the distributional sense) to the Dirac mass $\delta_0$ at zero as $t \rightarrow \infty$, then the  strong (Kelvin-Voigt) damping terms are recovered and system \eqref{MM-termic-full} formally  falls upon \eqref{BGN-termic-expanded}. Moreover, as suggested by \cite{gurtin-pipkin}, a formal integration by parts results
$$
\int_0^\infty \mu_i(s) p^i_{txx}(t-s) ds = \mu_i(0)p^i_{xx}(t) + \int_0^\infty \mu'_i(s) p^i_{xx}(t-s) ds,
$$
and regarding the notations $g_i:=-\mu_i'$  and $\mu_i=\beta +\varpi\mu_i(0)$, we can rewrite the last three equations  of
\eqref{MM-termic-full} as 
\begin{equation}\label{MM-heat-thermdisplacements}
\left\{\begin{array}{l}
	\varrho_1 p^1_{tt}  -\mu_1 \gamma_1 p^1_{xx}
+\varpi \gamma_1\displaystyle\int_0^\infty g_1(s) p^1_{xx}(t-s) ds + m_1(\varphi_x + \psi + \ell w)_t \, = \, 0,   \medskip\\
\varrho_2 p^2_{tt}  -\mu_2 \gamma_2 p^2_{xx}+\varpi \gamma_2 \displaystyle\int_0^\infty g_2(s) p^2_{xx}(t-s) ds+ m_2(w_x - \ell \varphi)_t \, = \, 0, \medskip \\ 
\varrho_3 p^3_{tt}  -\mu_3 \gamma_3 p^3_{xx}
+\varpi \gamma_3 \displaystyle\int_0^\infty g_3(s) p^3_{xx}(t-s) ds+ m_3\psi_{xt} \, = \, 0.
\end{array}\right.	
\end{equation}
Again, by means of \eqref{MM-law}, system \eqref{MM-heat-thermdisplacements} can be described in terms of the temperature as follows
\begin{equation}\label{MM-heatequations}
	\left\{	\begin{array}{l}
	\varrho_1 \theta^1_{tt}  -\mu_1 \gamma_1 \theta^1_{xx}
+\varpi \gamma_1\displaystyle\int_0^\infty g_1(s) \theta^1_{xx}(t-s) ds + m_1(\varphi_x + \psi + \ell w)_{tt} \, = \, 0,   \medskip\\
\varrho_2 \theta^2_{tt}  -\mu_2 \gamma_2 \theta^2_{xx}+\varpi \gamma_2 \displaystyle\int_0^\infty g_2(s) \theta^2_{xx}(t-s) ds+ m_2(w_x - \ell \varphi)_{tt} \, = \, 0, \medskip \\ 
\varrho_3 \theta^3_{tt}  -\mu_3 \gamma_3 \theta^3_{xx}
+\varpi \gamma_3 \displaystyle\int_0^\infty g_3(s) \theta^3_{xx}(t-s) ds+ m_3\psi_{xtt} \, = \, 0.
	\end{array}
	\right.	
\end{equation}

Finally, the thermoelastic Bresse-Type III system with memory relaxation is given by equations \eqref{MM-termic-full}$_{1,2,3}$ and \eqref{MM-heatequations}. Such a model (and the following ones) constitutes a new class of thermoelastic Bresse systems never considered in the literature. The remaining cases are similar and, for this reason, the details of computations will be omitted. 

\smallskip 

\noindent\textbf{II. BMR-Type III with double thermal coupling.}
\begin{equation}\label{MM-bending-axial}
\hskip-1.8cm \left\{	\begin{array}{l}
		\rho_1\varphi_{tt}  - k(\varphi_{x}+\psi + \ell w)_{x} - \ell k_0(w_x - \ell \varphi) + \ell m_2 \theta^2 \, = \, 0,  \medskip \\
		\rho_2\psi_{tt}  -  b\,\psi_{xx} +k(\varphi_{x}+\psi + \ell w) +m_3\theta^3_x  \, = \, 0, \medskip \\
		\rho_1 w_{tt} - k_0(w_x - \ell \varphi)_x + \ell k(\varphi_{x}+\psi + \ell w) +m_2\theta^2_x \,= \,0,  \medskip \\
	\varrho_2 \theta^2_{tt}  -\mu_2 \gamma_2 \theta^2_{xx}+\varpi \gamma_2 \displaystyle\int_0^\infty g_2(s) \theta^2_{xx}(t-s) ds+ m_2(w_x - \ell \varphi)_{tt} \, = \, 0, \medskip \\ 
	\varrho_3 \theta^3_{tt}  -\mu_3 \gamma_3 \theta^3_{xx}
	+\varpi \gamma_3 \displaystyle\int_0^\infty g_3(s) \theta^3_{xx}(t-s) ds+ m_3\psi_{xtt} \, = \, 0.
	\end{array}
	\right.	
\end{equation} 
\begin{equation}\label{MM-shear-axial}
	\left\{	\begin{array}{l}
		\rho_1\varphi_{tt}  - k(\varphi_{x}+\psi + \ell w)_{x} - \ell k_0(w_x - \ell \varphi) + m_1 \theta^1_x + \ell m_2 \theta^2 \, = \, 0,  \medskip \\
		\rho_2\psi_{tt}  -  b\,\psi_{xx} +k(\varphi_{x}+\psi + \ell w) - m_1 \theta^1  \, = \, 0, \medskip \\
		\rho_1 w_{tt} - k_0(w_x - \ell \varphi)_x + \ell k(\varphi_{x}+\psi + \ell w) +m_2\theta^2_x- \ell m_1\theta^1 \,= \,0,  \medskip \\
		\varrho_1 \theta^1_{tt}  -\mu_1 \gamma_1 \theta^1_{xx}
	+\varpi \gamma_1\displaystyle\int_0^\infty g_1(s) \theta^1_{xx}(t-s) ds + m_1(\varphi_x + \psi + \ell w)_{tt} \, = \, 0,   \medskip\\
	\varrho_2 \theta^2_{tt}  -\mu_2 \gamma_2 \theta^2_{xx}+\varpi \gamma_2 \displaystyle\int_0^\infty g_2(s) \theta^2_{xx}(t-s) ds+ m_2(w_x - \ell \varphi)_{tt} \, = \, 0.
	\end{array}
	\right.	
\end{equation} 
\begin{equation}\label{MM-shear-bending}
	\left\{	\begin{array}{l}
		\rho_1\varphi_{tt}  - k(\varphi_{x}+\psi + \ell w)_{x} - \ell k_0(w_x - \ell \varphi) + m_1 \theta^1_x \, = \, 0,  \medskip \\
		\rho_2\psi_{tt}  -  b\,\psi_{xx} +k(\varphi_{x}+\psi + \ell w) +m_3\theta^3_x - m_1 \theta^1  \, = \, 0, \medskip \\
		\rho_1 w_{tt} - k_0(w_x - \ell \varphi)_x + \ell k(\varphi_{x}+\psi + \ell w) - \ell m_1\theta^1 \,= \,0,  \medskip \\
	\varrho_1 \theta^1_{tt}  -\mu_1 \gamma_1 \theta^1_{xx}
+\varpi \gamma_1\displaystyle\int_0^\infty g_1(s) \theta^1_{xx}(t-s) ds + m_1(\varphi_x + \psi + \ell w)_{tt} \, = \, 0,   \medskip\\
\varrho_3 \theta^3_{tt}  -\mu_3 \gamma_3 \theta^3_{xx}
+\varpi \gamma_3 \displaystyle\int_0^\infty g_3(s) \theta^3_{xx}(t-s) ds+ m_3\psi_{xtt} \, = \, 0.
	\end{array}
	\right.	
\end{equation}

\noindent\textbf{III. BMR-Type III  with single thermal coupling.}
\begin{equation}\label{MM-shear}
	\left\{	\begin{array}{l}
		\rho_1\varphi_{tt}  - k(\varphi_{x}+\psi + \ell w)_{x} - \ell k_0(w_x - \ell \varphi) + m_1 \theta^1_x \, = \, 0,  \medskip \\
		\rho_2\psi_{tt}  -  b\,\psi_{xx} +k(\varphi_{x}+\psi + \ell w)  - m_1 \theta^1  \, = \, 0, \medskip \\
		\rho_1 w_{tt} - k_0(w_x - \ell \varphi)_x + \ell k(\varphi_{x}+\psi + \ell w) - \ell m_1\theta^1 \,= \,0,  \medskip \\
		\varrho_1 \theta^1_{tt}  -\mu_1 \gamma_1 \theta^1_{xx}
	+\varpi \gamma_1\displaystyle\int_0^\infty g_1(s) \theta^1_{xx}(t-s) ds + m_1(\varphi_x + \psi + \ell w)_{tt} \, = \, 0.
	\end{array}
	\right.	
\end{equation} 
\begin{equation}\label{MM-bending}
\hskip-3.15cm	\left\{	\begin{array}{l}
		\rho_1\varphi_{tt}  - k(\varphi_{x}+\psi + \ell w)_{x} - \ell k_0(w_x - \ell \varphi)  \, = \, 0,  \medskip \\
		\rho_2\psi_{tt}  -  b\,\psi_{xx} +k(\varphi_{x}+\psi + \ell w) +m_3\theta^3_x  \, = \, 0, \medskip \\
		\rho_1 w_{tt} - k_0(w_x - \ell \varphi)_x + \ell k(\varphi_{x}+\psi + \ell w) \,= \,0,  \medskip \\
	\varrho_3 \theta^3_{tt}  -\mu_3 \gamma_3 \theta^3_{xx}
	+\varpi \gamma_3 \displaystyle\int_0^\infty g_3(s) \theta^3_{xx}(t-s) ds+ m_3\psi_{xtt} \, = \, 0.
	\end{array}
	\right.	
\end{equation} 
\begin{equation}\label{MM-axial}
\hskip-1.85cm	\left\{	\begin{array}{l}
		\rho_1\varphi_{tt}  - k(\varphi_{x}+\psi + \ell w)_{x} - \ell k_0(w_x - \ell \varphi) + \ell m_2 \theta^2 \, = \, 0,  \medskip \\
		\rho_2\psi_{tt}  -  b\,\psi_{xx} +k(\varphi_{x}+\psi + \ell w)   \, = \, 0, \medskip \\
		\rho_1 w_{tt} - k_0(w_x - \ell \varphi)_x + \ell k(\varphi_{x}+\psi + \ell w) +m_2\theta^2_x \,= \,0,  \medskip \\
	\varrho_2 \theta^2_{tt}  -\mu_2 \gamma_2 \theta^2_{xx}+\varpi \gamma_2 \displaystyle\int_0^\infty g_2(s) \theta^2_{xx}(t-s) ds+ m_2(w_x - \ell \varphi)_{tt} \, = \, 0.
	\end{array}
	\right.	
\end{equation}

As in the previous case, the BMR-Type III models \eqref{MM-termic-full}-\eqref{MM-axial} constitute a new  class of thermoelastic Bresse systems, and yet not modeled in the previous literature. Moreover, the particular instance \eqref{MM-bending} with null history ($\theta^3\big|_{s<0}=0$) gives precisely the physical construction for the thermoelastic Bresse system made up recently by \cite{djellali-labidi-taallah}.

%
%

%

%
%
%
%

\paragraph{Funding.}
  M. A. Jorge has been  supported by the CNPq, grant \#160954/2022-3. T. F. Ma  has been   supported by the CNPq grant \#315165/2021-9 and FAPDF grant \#193.00001821/2022-21.

 \paragraph{Acknowledgments.}
The first  author  is very grateful to  the Professor To Fu Ma for
all kindness during his stay at the University of Bras\'ilia along the year 2023 and for all mathematical insights that led to the opportunity of developing the present scientific research.

\end{document}